\documentclass[11pt,letterpaper]{amsart}
\usepackage[T1]{fontenc}
\usepackage{lmodern}
\usepackage{amsmath,amssymb,amsthm}
\usepackage{enumitem}
\usepackage{microtype}
\usepackage{xcolor}

\theoremstyle{plain}
\newtheorem{theorem}{Theorem}[section]
\newtheorem{proposition}[theorem]{Proposition}
\newtheorem{lemma}[theorem]{Lemma}
\newtheorem{corollary}[theorem]{Corollary}
\newtheorem{conjecture}[theorem]{Conjecture}
\theoremstyle{definition}
\newtheorem{definition}[theorem]{Definition}

\theoremstyle{remark}

\DeclareMathOperator{\tr}{tr}
\DeclareMathOperator{\NS}{NS}
\DeclareMathOperator{\Pic}{Pic}
\DeclareMathOperator{\rank}{rank}
\DeclareMathOperator{\Def}{Def}
\DeclareMathOperator{\Psef}{Psef}
\DeclareMathOperator{\Bigcone}{Big}
\DeclareMathOperator{\Hom}{Hom}

\newcommand{\C}{\mathbb{C}}
\newcommand{\R}{\mathbb{R}}
\newcommand{\Z}{\mathbb{Z}}
\newcommand{\PP}{\mathbb{P}}
\newcommand{\dbar}{\bar\partial}
\newcommand{\pdbar}{\partial\dbar}
\newcommand{\cC}{\mathcal{C}}
\newcommand{\cE}{\mathcal{E}}

\newcommand{\cO}{\mathcal{O}}
\newcommand{\cK}{\mathcal{K}}
\newcommand{\cM}{\mathcal{M}}
\newcommand{\cX}{\mathcal{X}}
\newcommand{\cY}{\mathcal{Y}}
\newcommand{\cZ}{\mathcal{Z}}
\newcommand{\tX}{\tilde X}
\newcommand{\tpi}{\tilde\pi}

\begin{document}

\title[Smooth limits of compact K\"ahler manifolds with $c_1=0$]
{Smooth limits of compact K\"ahler manifolds\\ with vanishing first Chern class, I}

\author{Mu-Lin Li}
\address{School of Mathematics, Hunan University, Changsha 410082, China}
\email{mulin@hnu.edu.cn}


\begin{abstract}
Let $\pi\colon X\to\Delta$ be a proper holomorphic submersion from a smooth
complex manifold onto the unit disc. Suppose that for
every $t\in\Delta^*$ the fiber $X_t$ is a compact K\"ahler manifold whose real
first Chern class vanishes. We prove that the central fiber $X_0$ belongs to
Fujiki's class $\cC$, that is, $X_0$ is bimeromorphic to a compact K\"ahler
manifold.
\end{abstract}

\maketitle

\section{Introduction}

Let $\Delta\subset\C$ be a disc centered at the origin and let $\pi\colon X\longrightarrow\Delta$ be a proper holomorphic submersion. Denote by $
\Delta^*=\Delta\setminus\{0\}$ the punctured disc. The geometry of the
central fiber $X_0$ is constrained by that of the fibers over the punctured
disc, but the K\"ahler property itself need not survive smooth specialization:
an example of Hironaka \cite{Hironaka62} exhibits a smooth family of projective
threefolds whose central fiber is not K\"ahler. The general expectation is
nevertheless that a smooth limit of compact K\"ahler manifolds stays
bimeromorphic to a compact K\"ahler manifold, that is, belongs to the class
introduced by Fujiki \cite{Fujiki78}.

A compact complex manifold $Y$ belongs to \emph{Fujiki's class
$\cC$} if there exists a proper bimeromorphic map $Y'\dashrightarrow Y$ with
$Y'$ compact K\"ahler, or equivalently,  if $Y$
carries a K\"ahler current, i.e.\ a closed positive $(1,1)$-current which
dominates a Hermitian form. Class $\cC$ is stable under proper bimeromorphic
maps, under modifications, and under finite quotients; it is the natural
bimeromorphic closure of the K\"ahler condition.  This motivates the
following conjecture.

\begin{conjecture}\label{conj:intro}
Let $\pi\colon X\longrightarrow\Delta$ be a proper holomorphic submersion. If
$X_t$ is K\"ahler for every $t\in\Delta^*$, then the central fiber $X_0$
belongs to Fujiki's class $\cC$.
\end{conjecture}

Siu \cite{Siu83} proved that every K3 surface is K\"ahler, and hence that a
smooth limit of K3 surfaces is K\"ahler; see also Lamari \cite{Lamari99} for an
approach through currents. The author--Liu \cite{LL24}
proved that a limit of K\"ahler manifolds is K\"ahler under mild metric and
topological hypotheses on the central fiber, and the author
\cite{Li24} proved Conjecture~\ref{conj:intro} under the assumptions that $\pi$
is K\"ahler at one point of the punctured disc and that the upper mass of the
central fiber is finite.

The case in which the general fibers are projective or Moishezon has been
studied extensively. Popovici \cite{Popovici13} proved that if $X_t$ is
projective for every $t\in\Delta^*$ and the central fiber either has locally
deformation invariant $h^{0,1}$ or carries a strongly Gauduchon metric, then
$X_0$ is Moishezon. Barlet \cite{Barlet15} obtained analogous semicontinuity
results for the algebraic dimension of Moishezon fibers, and Rao--Tsai
\cite{RT21} proved related results under the existence of uncountably many
Moishezon fibers in the family. More recently, Rao \cite{Rao24} proved that a
limit of Moishezon manifolds is again Moishezon, without any further
hypothesis.

In this paper we treat the case in which the fibers over the punctured disc are
compact K\"ahler manifolds with vanishing real first Chern class. This is the
case of maximal relevance for the classification theory of compact
K\"ahler manifolds with $c_1=0$ (Calabi--Yau, hyperk\"ahler, and torus
manifolds), and it is the case in which the Beauville--Bogomolov decomposition
provides a substitute for the missing positivity.

\begin{theorem}\label{thm:main}
Let $\pi\colon X\to\Delta$ be a proper holomorphic submersion with connected
fibers $X_t$. Suppose that for every $t\in\Delta^*$ the fiber $X_t$ is
K\"ahler and $c_1(X_t)=0$ in $H^2(X_t,\R)$. Then $X_0$ belongs to Fujiki's
class $\cC$.
\end{theorem}

Theorem~\ref{thm:main} settles, in the case of vanishing first Chern class, the
smooth one-parameter specialization question raised by Soldatenkov and
Verbitsky \cite[Conjecture 1.2]{SV24}, and it does so without assuming the
existence of a holomorphic symplectic form on the central fiber. Liu and Shen
\cite[Theorem 2]{LS26} prove a bimeromorphic hyperk\"ahler limit theorem under
the additional hypothesis that the central fiber carries a strongly Gauduchon
metric.

We sketch here our method as following. Throughout, we work on a finite
\'etale cover of $X$, which is harmless for the conclusion because class $\cC$
is stable under finite \'etale maps in both directions.

First, we prove the relative Beauville--Bogomolov product. On a suitable finite
\'etale cover $\tX\to X$ the Beauville--Bogomolov factors of the punctured
fibers can be assembled into global families over $\Delta^*$
(Theorem~\ref{thm:productcover}), and the resulting relative product extends
across the central parameter (Theorem~\ref{thm:product}). The extension is
obtained by realizing the factors as leaf spaces of holomorphic foliations
defined by the limiting factor forms; boundedness of the leaf volumes, which
follows from the smoothness of the extended cohomology classes, makes Bishop's
compactness theorem available and yields properness of the foliations.

Then we study the extension of the periods. The top degree period of each factor
extends as a nonzero class across the central parameter, and the same holds for
the pairing of that class with its conjugate, because the relative canonical
class is trivial and the product of all factors degenerates only to finite
order. This produces a smooth extension of the degree two Hodge splitting
inside the flat cohomology bundle, hence a smooth family of classes which is
K\"ahler on a dense open subset of the punctured disc and has constant volume.

A version of the mass estimate of
\cite{Li24,Rao24} which is uniform in the parameter, together with Bochner's
principle for Ricci-flat K\"ahler metrics, shows that the normalized factor
forms have no poles along $X_0$: the mass of a K\"ahler representative of the
moving class is bounded, and a genuine pole would force the integral of a
suitable power of the pointwise norm of the factor form to blow up.

Finally, we study the three types of factors. A torus factor extends as a torus. For
a simply connected factor with vanishing intermediate holomorphic forms, the
Moishezon property of the central fiber follows by extending a fixed integral
K\"ahler class to a line bundle and applying Riemann--Roch; this is a direct
argument, and Rao's general specialization theorem \cite[Corollary
1.4]{Rao24} applies as well. For an irreducible holomorphic symplectic factor
we perturb the nearby fibers to Picard number zero, use the surjectivity of the
period map on a connected component of the marked moduli space to compare them
with a genuine smooth K\"ahler family, and deduce from a meromorphic
compactness theorem for holomorphic maps of bounded energy that the central
fiber is bimeromorphic to a hyperk\"ahler manifold.

\section{Preliminaries}\label{sec:prelim}

\subsection{The Beauville--Bogomolov decomposition}

In this paper a \emph{Calabi--Yau manifold} is a compact K\"ahler manifold $X$
with $K_X\cong\cO_X$; a \emph{strict Calabi--Yau manifold} is a simply
connected Calabi--Yau manifold $X$ of dimension at least three with $H^0(X,\Omega^p_X)=0$ for
$0<p<\dim X$. We also recall the hyperk\"ahler terminology and the period
theory used below; our conventions follow \cite{Beauville83,Fujiki87,Huy12}.

\begin{definition}\label{def:ihs}
An \emph{irreducible compact hyperk\"ahler manifold}, also called an
\emph{irreducible holomorphic symplectic manifold}, is a simply connected
compact K\"ahler manifold $Z$ such that $H^0(Z,\Omega^2_Z)=\C\sigma$ for a
two-form $\sigma$ which is everywhere nondegenerate.
\end{definition}

The complex dimension of $Z$ is even, say $d=2m$. Holomorphic forms on a
compact K\"ahler manifold are closed, so $d\sigma=0$, and the nowhere vanishing
form $\sigma^m$ trivializes $K_Z$. Simple connectedness together with Hodge
theory gives $H^1(Z,\cO_Z)=0$ and $h^{2,0}(Z)=1$.

Contraction with $\sigma$ defines an isomorphism of holomorphic vector bundles
\begin{equation}\label{eq:contraction}
T_Z\xrightarrow{\ \sim\ }\Omega^1_Z,\qquad v\longmapsto\iota_v\sigma .
\end{equation}
In particular $H^0(Z,T_Z)=0$ and $H^1(Z,T_Z)\cong H^1(Z,\Omega^1_Z)$.

By Yau's theorem \cite{Yau78} every K\"ahler class on $Z$ contains a unique
Ricci-flat K\"ahler form. If $\hat\eta$ is such a form, then by Bochner's
principle every holomorphic $p$-form on $Z$ is parallel with respect to the
associated K\"ahler metric \cite{Beauville83}; in particular $\sigma$ is
parallel, so the restricted holonomy of the metric is contained in
$\mathrm{Sp}(m)$. Moreover the holonomy representation is irreducible: a
holonomy invariant splitting of $T_Z$ would decompose $\sigma$ as a sum of
parallel two-forms which, being of pure bidegree, are holomorphic, and this
contradicts $h^{2,0}(Z)=1$ together with the nondegeneracy of $\sigma$. Thus
the metric is de Rham irreducible.

The following structure theorem is the starting point of everything that
follows.

\begin{theorem}[{\cite[Theorem 2]{Beauville83}}]\label{thm:BB}
Let $X$ be a compact K\"ahler manifold with $c_1(X)=0$ in $H^2(X,\R)$. Then
there exist a finite \'etale cover $p\colon\tX\to X$ and an isomorphism
$\tX\cong T\times\prod_i X_i\times\prod_j Y_j$ such that:
\begin{enumerate}[label=\textup{(\roman*)}]
\item $T$ is a complex torus;
\item $X_i$ is a strict Calabi--Yau manifold of dimension at least $3$, that
is, $K_{X_i}$ is trivial and $H^0(X_i,\Omega^p_{X_i})=0$ for all
$0<p<\dim X_i$;
\item $Y_j$ is an irreducible compact hyperk\"ahler manifold.
\end{enumerate}
\end{theorem}

We shall refer to $T$, the $X_i$ and the $Y_j$ as the
\emph{Beauville--Bogomolov factors} of $X$, and to the decomposition itself as
the \emph{Beauville--Bogomolov decomposition}. Notice that all factors are
simply connected except possibly the torus, and that at most one torus occurs.

\subsection{The Beauville--Bogomolov form and the period domain}

Let $Z$ be an irreducible compact hyperk\"ahler manifold of dimension $d=2m$.
Since $Z$ is simply connected we have $H_1(Z,\Z)=0$, so the universal
coefficient theorem shows that $H^2(Z,\Z)$ is torsion free; it carries the
primitive integral \emph{Beauville--Bogomolov quadratic form} $q_Z$, with sign
chosen so that $q_Z$ is positive on K\"ahler classes. We write
$q_Z(\alpha,\beta)$ for the associated symmetric bilinear form and
$q_Z(\alpha)=q_Z(\alpha,\alpha)$. Its real signature is $(3,b_2(Z)-3)$, and
there is a positive real number $c_Z$, the \emph{Fujiki constant}, such that
\begin{equation}\label{eq:fujiki}
\int_Z\alpha^{2m}=c_Z\,q_Z(\alpha)^m
\qquad(\alpha\in H^2(Z,\C)).
\end{equation}
These statements are due to Beauville \cite{Beauville83} and Fujiki
\cite{Fujiki87}. Parallel transport in a family of irreducible compact
hyperk\"ahler manifolds preserves both $q_Z$ and $c_Z$, so both are invariants
of the deformation type.

We shall use the following standard consequences; see
\cite[Section 1.2]{Huy12} for a unified treatment.

\begin{proposition}\label{prop:qproperties}
Let $Z$ be an irreducible compact hyperk\"ahler manifold of dimension $d=2m$.
\begin{enumerate}[label=\textup{(\roman*)}]
\item The Hodge decomposition is orthogonal for $q_Z$ in the following sense:
\begin{equation}\label{eq:qperiod}
q_Z([\sigma])=0,\qquad q_Z([\sigma],[\bar\sigma])>0,
\end{equation}
\begin{equation}\label{eq:qHodge}
H^{1,1}(Z,\C)=\bigl(\C[\sigma]\oplus\C[\bar\sigma]\bigr)^{\perp_{q_Z}} .
\end{equation}
\item The real plane spanned by $\operatorname{Re}[\sigma]$ and
$\operatorname{Im}[\sigma]$ is positive definite, and the restriction of $q_Z$
to $H^{1,1}(Z,\R)$ has signature $(1,b_2(Z)-3)$.
\item For classes $x,y\in H^2(Z,\C)$ with $q_Z(x)=q_Z(y)=0$ one has
\begin{equation}\label{eq:polarizedfujiki-main}
\binom{2m}{m}\int_Z x^my^m
=c_Z\bigl(2q_Z(x,y)\bigr)^{m},
\end{equation}
and more generally the integral of any monomial of degree $2m$ in classes of
degree two is a polynomial, with coefficients depending only on $m$, in the
pairings $q_Z$ of those classes.
\end{enumerate}
\end{proposition}

Assertion \textup{(iii)} is the polarized form of \eqref{eq:fujiki}: it follows
by applying \eqref{eq:fujiki} to $x+sy$ and comparing the coefficients of
$s^m$. The last sentence is the content of Fujiki's relations \cite{Fujiki87}.

Fix a lattice $(\Lambda,q_\Lambda)$ isometric to $(H^2(Z,\Z),q_Z)$ and call a
lattice isometry $\phi\colon H^2(Z,\Z)\xrightarrow{\sim}\Lambda$ a
\emph{marking}. Two marked manifolds $(Z,\phi)$ and $(Z',\phi')$ are
\emph{isomorphic} if there is a biholomorphism $f\colon Z\to Z'$ with
$\phi\circ f^*=\phi'$. Put $\Lambda_\C=\Lambda\otimes_\Z\C$ and
$\Lambda_\R=\Lambda\otimes_\Z\R$. The \emph{period domain} is the open subset
\begin{equation}\label{eq:perioddomain}
D_\Lambda=\bigl\{[\ell]\in\PP(\Lambda_\C) :
q_\Lambda(\ell)=0,\ q_\Lambda(\ell,\bar\ell)>0\bigr\}
\end{equation}
of the nonsingular quadric $\{q_\Lambda=0\}\subset\PP(\Lambda_\C)$; it is a
homogeneous domain of complex dimension $b_2(Z)-2$. The \emph{period} of
$(Z,\phi)$ is the point $[\phi_\C([\sigma])]$. It records the line
$H^{2,0}(Z)$, whose complex conjugate gives $H^{0,2}(Z)$, while the
$q$-orthogonal complement of the sum of these two lines is $H^{1,1}(Z)$, by
\eqref{eq:qHodge}.

Over a contractible base a marking at one fiber extends by Gauss--Manin
parallel transport. Once all cohomology groups have been identified with the
same lattice we suppress the marking from the notation and write $q$ for the
transported form.

\subsection{Local deformations and the local Torelli theorem}

Let $\cZ\to\Def(Z)$ be a sufficiently small Kuranishi family with distinguished
fiber $Z$. Its tangent space at the distinguished point is $H^1(Z,T_Z)$, and
the vanishing $H^0(Z,T_Z)=0$ ensures that the family may be taken locally
universal. After extending a marking by parallel transport, the line of
holomorphic two-forms defines the holomorphic \emph{period map}
\[
P\colon\Def(Z)\longrightarrow D_\Lambda,\qquad
s\longmapsto\bigl[\phi_{s,\C}([\sigma_s])\bigr].
\]

\begin{theorem}\label{thm:localTorelli}
The germ $\Def(Z)$ is smooth of dimension $h^{1,1}(Z)=b_2(Z)-2$, and its
sufficiently small fibers are irreducible compact hyperk\"ahler manifolds.
After shrinking the base, $P$ is biholomorphic onto an open neighborhood of the
period of $Z$ in $D_\Lambda$.
\end{theorem}

\begin{proof}
Unobstructedness follows from the Bogomolov--Tian--Todorov theorem
\cite{Tian87,Todorov89}; the dimension statement follows from
\eqref{eq:contraction}, which gives
$h^1(Z,T_Z)=h^1(Z,\Omega^1_Z)=h^{1,1}(Z)=b_2(Z)-2$. We recall the differential
argument for the local Torelli statement, referring to \cite{Beauville83} and
\cite[Section 4.1 and Theorem 4.1]{Huy12} for details.

Using the marking and the Hodge decomposition,
\[
T_{[\sigma]}D_\Lambda
=\Hom\bigl(\C[\sigma],[\sigma]^{\perp_q}/\C[\sigma]\bigr)
\ \cong\ \Hom\bigl(H^{2,0}(Z),H^{1,1}(Z)\bigr),
\]
where we used \eqref{eq:qHodge} and the fact that $[\sigma]^{\perp_q}$ is
spanned by $H^{1,1}(Z)$ together with $[\sigma]$. With the contraction
convention for the Kodaira--Spencer map, the differential is
\begin{equation}\label{eq:dP}
dP_0(\xi)\bigl([\sigma]\bigr)=[\xi\lrcorner\,\sigma],
\qquad \xi\in H^1(Z,T_Z).
\end{equation}
By \eqref{eq:contraction}, evaluation at $[\sigma]$ identifies
\eqref{eq:dP} with the isomorphism
$H^1(Z,T_Z)\xrightarrow{\sim}H^1(Z,\Omega^1_Z)\cong H^{1,1}(Z)$. The inverse
function theorem now gives the assertion about $P$. Stability of the
irreducible compact hyperk\"ahler property under small deformations is part of
the same classical theory.
\end{proof}

\subsection{Global marked moduli and componentwise periods}

Fix the deformation type of a reference irreducible compact hyperk\"ahler
manifold $Z^*$ and a lattice $(\Lambda,q_\Lambda)\cong(H^2(Z^*,\Z),q_{Z^*})$.
The \emph{full marked moduli space} of this deformation type is
\[
\cM_\Lambda=\{(Z,\phi)\}/\!\sim,
\]
where $Z$ runs over irreducible compact hyperk\"ahler manifolds deformation
equivalent to $Z^*$ and $\phi$ over all lattice isometries
$H^2(Z,\Z)\xrightarrow{\sim}\Lambda$, and where $(Z,\phi)\sim(Z',\phi')$ if
there is a biholomorphism $f\colon Z\to Z'$ with $\phi\circ f^*=\phi'$. The
deformation type is fixed once and for all and is suppressed from the notation
$\cM_\Lambda$. Marked local deformation spaces provide complex manifold charts
for $\cM_\Lambda$; this space may fail to be Hausdorff and may have more than
one connected component. We write $\cM_\Lambda=\coprod_\nu\cM^\nu_\Lambda$ and
call a \emph{connected marked deformation component} one such component. The
Hodge bundles are smooth bundles in each local chart, and markings extend by
Gauss--Manin parallel transport.

The \emph{global period map} is
\[
\mathcal P\colon\cM_\Lambda\longrightarrow D_\Lambda,\qquad
(Z,\phi)\longmapsto\PP\bigl(\phi_\C(H^{2,0}(Z))\bigr).
\]
It is well defined on marked isomorphism classes, and its restriction to a
local deformation chart is the local period map of
Theorem~\ref{thm:localTorelli}. For each connected component we put
$\mathcal P_\nu:=\mathcal P|_{\cM^\nu_\Lambda}\colon
\cM^\nu_\Lambda\to D_\Lambda$.

\begin{theorem}[{\cite[Theorems 4.1, 5.5 and 1.3]{Huy12}}]
\label{thm:globalTorelli}
For the fixed deformation type above the following hold.
\begin{enumerate}[label=\textup{(\roman*)}]
\item The map $\mathcal P$ is locally biholomorphic.
\item For every connected component $\cM^\nu_\Lambda$ the map
$\mathcal P_\nu$ is surjective onto the whole period domain $D_\Lambda$;
consequently $\mathcal P$ itself is surjective.
\item Let
\[
D^{\mathrm{gen}}_\Lambda
:=D_\Lambda\setminus\bigcup_{0\neq\lambda\in\Lambda}
\bigl\{[\ell]\in D_\Lambda : q_\Lambda(\lambda,\ell)=0\bigr\}.
\]
For every $p\in D^{\mathrm{gen}}_\Lambda$ and every fixed connected component,
the fiber $\mathcal P_\nu^{-1}(p)$ consists of exactly one marked isomorphism
class.
\end{enumerate}
\end{theorem}

Thus every period $p\in D_\Lambda$, including a period lying on an integral
Hodge hyperplane, is realized by a marked hyperk\"ahler manifold in any
prescribed component. The singleton-fiber assertion is made within a fixed
component and at the specified generic periods; by the Hodge theoretic
description below these are precisely the periods of Picard number zero
manifolds. In particular, if $(Z,\phi)$ and $(Y,\psi)$ lie in the same
component, have Picard number zero, and have the same period, then there is a
biholomorphism $f\colon Z\to Y$ satisfying $\phi\circ f^*=\psi$.

\subsection{Integral Hodge loci and positive cones}

The N\'eron--Severi group and the Picard number of an irreducible compact
hyperk\"ahler manifold are
\[
\NS(Z)=H^2(Z,\Z)\cap H^{1,1}(Z),\qquad
\rho(Z)=\rank_\Z\NS(Z).
\]
The exponential sequence together with the Lefschetz $(1,1)$ theorem
identifies $\Pic(Z)$ with $\NS(Z)$, since $H^1(Z,\cO_Z)=0$. For a nonzero class
$\lambda\in\Lambda$ the \emph{period hyperplane} is
\[
H_\lambda=\bigl\{[\ell]\in D_\Lambda : q_\Lambda(\lambda,\ell)=0\bigr\}.
\]
A real integral class is of type $(1,1)$ precisely when it is $q$-orthogonal
to the period: indeed, for real $\lambda$ the identity
$q_\Lambda(\lambda,\bar\ell)=\overline{q_\Lambda(\lambda,\ell)}$ shows that
$q_\Lambda(\lambda,\ell)=0$ is equivalent to $q_\Lambda(\lambda,\bar\ell)=0$,
and \eqref{eq:qHodge} then places $\lambda$ in $H^{1,1}$. Hence
\begin{equation}\label{eq:NSzero}
\rho(Z)=0\quad\Longleftrightarrow\quad
[\phi_\C([\sigma])]\notin\bigcup_{0\neq\lambda\in\Lambda}H_\lambda .
\end{equation}

\begin{corollary}\label{cor:density}
In a sufficiently small full local deformation space of an irreducible compact
hyperk\"ahler manifold the parameters with Picard number zero form the
complement of a countable union of closed proper analytic subsets. In
particular they meet every nonempty open subset of the local deformation space.
\end{corollary}

\begin{proof}
By Theorem~\ref{thm:localTorelli} it suffices to work in an open subset of
$D_\Lambda$. For $\lambda\neq0$ the linear functional
$q_\Lambda(\lambda,\cdot)$ is nonzero, and its projective hyperplane cannot
contain an open subset of the nondegenerate quadric $\{q_\Lambda=0\}$, which is
irreducible and is not contained in any hyperplane. Hence each $H_\lambda$ is a
closed proper analytic subset of the period domain, and its intersection with a
local period chart is nowhere dense. Since $\Lambda$ is countable, Baire's
theorem gives the density assertion.
\end{proof}

The \emph{K\"ahler cone} $\cK(Z)$ consists of the real $(1,1)$-classes
represented by K\"ahler forms. The \emph{positive cone} $\cC^+_Z$ is the
component of $\{a\in H^{1,1}(Z,\R):q_Z(a)>0\}$ containing $\cK(Z)$. The
Lorentzian signature shows that, for any fixed $\kappa\in\cK(Z)$,
\[
\cC^+_Z=\bigl\{a\in H^{1,1}(Z,\R):q_Z(a)>0,\ q_Z(a,\kappa)>0\bigr\}.
\]
A real $(1,1)$-class is \emph{pseudoeffective} if it contains a closed positive
$(1,1)$-current, and \emph{big} if it contains a K\"ahler current; we denote
these cones by $\Psef(Z)$ and $\Bigcone(Z)$. The two hyperk\"ahler cone
properties used below are
\begin{equation}\label{eq:cone}
\cC^+_Z\subset\Bigcone(Z)\subset\Psef(Z),
\qquad \rho(Z)=0\ \Longrightarrow\ \cK(Z)=\cC^+_Z .
\end{equation}
The first inclusion is \cite[Proposition 1]{Huy03a}, which rests on the
Demailly--P\u{a}un theorem; the second follows from hyperk\"ahler cone theory
\cite{Huy03b} and is also stated in \cite[Theorem 5.1 and Corollary 5.3]{Huy12}.
In particular, a class of positive square must lie in the correct positive
component before the Picard number zero conclusion can be used to make it
K\"ahler.

\section{Gauduchon masses and extension of cohomology classes}\label{sec:mass}

Throughout this section, a family is a proper holomorphic submersion with
connected fibers. Smooth families of differential forms are understood in a
differentiable trivialization.

\subsection{Gauduchon forms in families}

A Hermitian form $h$ on a compact complex $n$-fold is Gauduchon if
$\pdbar h^{n-1}=0$. Gauduchon's theorem \cite{Gauduchon77} gives such a form
in every conformal class. The following lemma for family of Gauduchon metrics is known for experts, we include a proof here for convenience.

\begin{lemma}\label{lem:parametricGauduchon}
Let $f\colon\cX\to S$ be a family of compact complex $n$-folds, and let $h_0$
be a Gauduchon form on $\cX_{s_0}$. After shrinking $S$ about $s_0$, there is
a smooth family of Gauduchon forms $h_s$ with $h_{s_0}=h_0$.
\end{lemma}

\begin{proof}
Choose a smooth family of Hermitian forms $\widehat h_s$ through $h_0$.
For $n\ge2$, consider the real elliptic operators
\[
 L_su=\frac{i\pdbar(u\widehat h_s^{n-1})}{\widehat h_s^{\,n}}.
\]
Their kernels have dimension one and are generated by positive functions,
by Gauduchon's theorem. The kernels therefore form a smooth line bundle:
one can, for example, use the spectral projections of $L_s^*L_s$ onto the
isolated eigenvalue zero, with adjoints taken in fixed smooth Sobolev
trivializations. Applying these projections to the constant function $1$
gives, after shrinking, a smooth positive family $u_s$ with $u_{s_0}=1$.
Elliptic regularity gives smoothness in both the fiber and the parameter.
Then $h_s=u_s^{1/(n-1)}\widehat h_s$ has the required properties.
\end{proof}

\subsection{A scalar phase estimate}

The distinction between boundedness and smooth extendibility in the next
lemma is essential.

\begin{lemma}\cite[Lemma 4.1]{Li24}\label{lem:scalarphase}
Let $A$ be a smooth real-valued function and $B$ a smooth complex-valued
function on a disc about $0$. Let $N\ge0$ be an integer. If
\[
 A(t)+2\operatorname{Re}\bigl(t^{-N}B(t)\bigr)\ge0
 \qquad(t\ne0),
\]
then $B(t)=O(|t|^N)$. In particular, the displayed function is bounded near
$0$ on the punctured disc.
\end{lemma}

\begin{proof}
There is nothing to prove when $N=0$. Otherwise, suppose that the Taylor
polynomial of $B$ of degree $N-1$ is nonzero, and let $k<N$ be its first
nonzero homogeneous degree. Write its homogeneous part as
$\sum_{p+q=k}b_{pq}t^p\bar t^q$. Along $t=re^{i\theta}$, the leading term
of $2\operatorname{Re}(t^{-N}B(t))$ is $2r^{k-N}\varphi(\theta)$, where
\[
 \varphi(\theta)=\operatorname{Re}
       \sum_{p+q=k}b_{pq}e^{i(p-q-N)\theta}.
\]
The frequencies $p-q-N$ are distinct and strictly negative. Consequently,
$\varphi$ is nonzero and has mean zero, and hence is negative at some
$\theta_0$. Along that ray, the displayed inequality fails for small $r$,
because $A$ is bounded and $k-N<0$. Thus all Taylor terms of degree less
than $N$ vanish, which gives the assertion by Taylor's theorem.
\end{proof}

\subsection{The mass of a pseudoeffective class}

For a Gauduchon form $h$ on a compact complex $n$-fold $Y$, integration
against $h^{n-1}$ defines a linear functional on real Bott--Chern
cohomology:
\[
 \mu_{Y,h}(b)=\int_Y T\wedge h^{n-1},
 \qquad b\in H^{1,1}_{\mathrm{BC}}(Y,\R).
\]
The value is independent of the representative, since a change by
$i\pdbar\psi$ integrates to zero. When $Y$ is K\"ahler, the
$\pdbar$-lemma identifies this group with the real de Rham classes of type
$(1,1)$, so the mass is then determined by the de Rham class.  A
pseudoeffective class has nonnegative mass.

\begin{proposition}\label{prop:mass}
Let $\pi\colon X\to\Delta$ be a family of compact complex $n$-folds whose
fibers over $\Delta^*$ are K\"ahler. Let $\beta(t)$ be a smooth real section
of the marked degree-two cohomology bundle over $\Delta$, of type $(1,1)$
over $\Delta^*$. Suppose that $\beta(t)$ is pseudoeffective on a dense
subset $U'\subset\Delta^*$. For every smooth family of Gauduchon forms
$h_t$ over $\Delta$, there is a nonnegative smooth function $m$ on a smaller
punctured disc, bounded near $0$, such that
\[
 \int_{X_t}T_t\wedge h_t^{n-1}=m(t)
\]
for every $t\in U'$ and every closed positive current $T_t\in\beta(t)$.
\end{proposition}

\begin{proof}
Fix a differentiable trivialization and smooth closed real two-forms
$\eta_1,\ldots,\eta_b$ on $X$ representing a flat basis of degree-two
cohomology. Write
\[
 \beta(t)=\sum_i c_i(t)[\eta_i|_{X_t}],
 \qquad \zeta_t=\sum_i c_i(t)\eta_i|_{X_t},
\]
where the $c_i$ are smooth real functions on $\Delta$.

Set $\mathcal F=R^2\pi_*\cO_X$ and
$Q=\mathcal F/\operatorname{Tor}(\mathcal F)$. On the smooth curve
$\Delta$, the coherent torsion-free sheaf $Q$ is locally free. The absolute
forms $\eta_i^{0,2}$ are $\dbar$-closed and define a holomorphic morphism
\[
 \Phi\colon\cO_\Delta^{\oplus b}\longrightarrow Q.
\]
Here we use the Stein--Leray identification
$H^2(X,\cO_X)=H^0(\Delta,R^2\pi_*\cO_X)$, and its local counterparts.
Over $\Delta^*$, the Hodge numbers are locally constant because all fibers
are K\"ahler; thus $\mathcal F$ is locally free there and satisfies base
change.

Let $\cE=\ker\Phi$. Its quotient in $\cO_\Delta^{\oplus b}$ is the
coherent subsheaf $\operatorname{im}\Phi$ of $Q$, hence is torsion free and
locally free. Therefore $\cE$ is a holomorphic subbundle. The vector
$c(t)=(c_i(t))$ lies in this subbundle for $t\ne0$, by the type assumption,
and hence also at $0$ by continuity. After shrinking, choose a holomorphic
frame $v_1,\ldots,v_s$ of $\cE$. There are smooth functions $\alpha_a$ on
$\Delta$ such that
\[
 c(t)=\sum_a\alpha_a(t)v_a(t).
\]
No holomorphicity is asserted for these coefficients.

For each $a$, the section of $\mathcal F$ defined by
$\sum_i v_{ai}(t)\eta_i^{0,2}$ is torsion supported at $0$. A common power
$t^N$ annihilates all these sections. Stein--Leray and the Dolbeault
resolution then give smooth absolute $(0,1)$-forms $b_a$ with
\[
 \dbar b_a=t^N\sum_i v_{ai}(t)\eta_i^{0,2}.
\]
Put $b=\sum_a\alpha_a b_a$ and $b_t=b|_{X_t}$. Although $\alpha_a$ need not
be holomorphic on the base, they are constant along each fiber, so
\[
 \dbar_{X_t}b_t=t^N\zeta_t^{0,2}.
\]
Consequently, for $t\ne0$, the form
\begin{equation}\label{eq:etaprime}
 \eta'_t=\zeta_t-d_{X_t}
      \bigl(t^{-N}b_t+\bar t^{-N}\overline{b_t}\bigr)
\end{equation}
is real, closed, of type $(1,1)$, and represents $\beta(t)$.

Its Gauduchon mass is
\begin{equation}\label{eq:massexp}
 m(t)=\int_{X_t}\eta'_t\wedge h_t^{n-1}
     =A(t)+2\operatorname{Re}\bigl(t^{-N}B(t)\bigr),
\end{equation}
where
\[
 A(t)=\int_{X_t}\zeta_t\wedge h_t^{n-1},\qquad
 B(t)=-\int_{X_t}\partial_{X_t}b_t\wedge h_t^{n-1}.
\]
Both $A$ and $B$ are smooth on $\Delta$, with $A$ real-valued. On each
K\"ahler fiber, the $\pdbar$-lemma for currents shows that this is the
mass of every positive representative of $\beta(t)$. Thus $m\ge0$ on
$U'$, and hence on $\Delta^*$ by continuity. Lemma~\ref{lem:scalarphase}
gives $B(t)=O(|t|^N)$ and proves the required boundedness.
\end{proof}

\subsection{Extension of a fixed class and the Moishezon case}

\begin{proposition}\label{prop:extendclass}
Let $\pi\colon X\to\Delta$ be a family whose fibers over $\Delta^*$ are
K\"ahler. Let $c\in H^2(X,\R)$ have restrictions of type $(1,1)$ and
pseudoeffective on a dense subset of $\Delta^*$. After shrinking $\Delta$,
the class $c$ has a smooth closed real $(1,1)$-representative on $X$. If
$c$ is the real image of $\widetilde c\in H^2(X,\Z)$, then
$\widetilde c=c_1(L)$ for a holomorphic line bundle $L$ on $X$.
\end{proposition}

\begin{proof}
Choose a smooth closed real representative $\eta$ of $c$, and let
$o(c)\in H^2(X,\cO_X)$ be the class represented by $\eta^{0,2}$. Under
Stein--Leray, it is a holomorphic section of $R^2\pi_*\cO_X$. Over
$\Delta^*$, this direct image is locally free and satisfies base change.
The section vanishes on the dense set where $c|_{X_t}$ is of type
$(1,1)$, hence vanishes on all of $\Delta^*$. It is therefore annihilated
by a power of $t$. Let $N\ge0$ be the least exponent with $t^No(c)=0$.

Suppose that $N>0$, and choose a smooth absolute $(0,1)$-form $b$ with
\[
 \dbar b=t^N\eta^{0,2}.
\]
The case $n=1$ cannot occur here, since then $R^2\pi_*\cO_X=0$; thus
$n\ge2$. Let $\Gamma_0$ be any strictly positive real
$(n-1,n-1)$-form on $X_0$ with $\pdbar\Gamma_0=0$. Such a form is the
$(n-1)$-st power of a unique Hermitian form. Lemma~\ref{lem:parametricGauduchon}
extends it to a smooth family $\Gamma_t=h_t^{n-1}$ with
$\pdbar\Gamma_t=0$.

The corrected representatives
\[
 \eta|_{X_t}-d_{X_t}
       \bigl(t^{-N}b_t+\bar t^{-N}\overline{b_t}\bigr),
 \qquad t\ne0,
\]
have nonnegative Gauduchon mass, first on the dense pseudoeffective locus
and then on all of $\Delta^*$. Applying Lemma~\ref{lem:scalarphase} to
this mass gives
\begin{equation}\label{eq:B0}
 \int_{X_0}\partial b_0\wedge\Gamma_0=0.
\end{equation}
Indeed, the corresponding smooth numerator is
$-\int_{X_t}\partial b_t\wedge\Gamma_t$, and its value at $0$ must vanish
because $N>0$.

Since $\dbar b_0=0$, the $(1,1)$-form $\partial b_0$ is $d$-closed.
Every real $\pdbar$-closed $(n-1,n-1)$-form is a difference of two strictly
positive such forms. Thus \eqref{eq:B0}, as $\Gamma_0$ varies, says that
the Bott--Chern class of $\partial b_0$ annihilates all complementary
Aeppli cohomology. The nondegenerate Bott--Chern--Aeppli pairing
\cite[Part~I, \S2.c]{Schweitzer07} yields a smooth complex-valued function
$u$ such that
\begin{equation}\label{eq:partialb0}
 \partial b_0=\pdbar u.
\end{equation}
Hence $\gamma=b_0-\dbar u$ is a $d$-closed complex $(0,1)$-form.

Ehresmann's theorem identifies $H^1(X,\C)$ with $H^1(X_0,\C)$, so the
de Rham class of $\gamma$ extends to $X$. Taking the $(0,1)$-part of a
smooth closed representative on $X$ shows that
$[b_0]=[\gamma]\in H^1(X_0,\cO_{X_0})$ extends to
$H^1(X,\cO_X)$.

Let $i\colon X_0\hookrightarrow X$. The connecting morphism in
\[
 0\longrightarrow\cO_X\xrightarrow{\ t\ }\cO_X
   \longrightarrow i_*\cO_{X_0}\longrightarrow0
\]
sends $[b_0]$ to $t^{N-1}o(c)$. To justify this computation with smooth
forms, take the quotient complex
$\mathcal Q^\bullet=A_X^{0,\bullet}/tA_X^{0,\bullet}$. Multiplication by
$t$ is injective termwise, and the long exact sequence of cohomology
sheaves shows that $\mathcal Q^\bullet$ is a fine resolution of
$i_*\cO_{X_0}$. Restriction induces a quasi-isomorphism
$\mathcal Q^\bullet\to i_*A_{X_0}^{0,\bullet}$, since it induces the
identity on the resolved sheaf. In particular, the quotient complex is
not being identified termwise with the complex of forms on $X_0$.
The class of $b$ in $\mathcal Q^1$ corresponds to $[b_0]$, and its
connecting class is represented by
$t^{-1}\dbar b=t^{N-1}\eta^{0,2}$, as claimed. Since $[b_0]$ extends to
$X$, exactness forces this class to vanish, contradicting minimality of
$N$.

It follows that $o(c)=0$. A smooth solution of
$\dbar b=\eta^{0,2}$ gives the closed real $(1,1)$-representative
$\eta-d(b+\bar b)$. Finally, the image of $\widetilde c$ in
$H^2(X,\cO_X)$ agrees with $o(c)$: the comparison map factors through
complex de Rham cohomology. The exponential sequence therefore realizes
$\widetilde c$ as the first Chern class of a holomorphic line bundle.
\end{proof}

\begin{corollary}\label{cor:strictCY}
Suppose that the fibers over $\Delta^*$ are compact K\"ahler manifolds
with $K_{X_t}\cong\cO_{X_t}$ and $h^{2,0}(X_t)=0$. Then $X_0$ is
Moishezon. In particular, this applies to families of strict Calabi--Yau
manifolds of dimension at least three.
\end{corollary}

\begin{proof}
The real degree-two cohomology of every punctured fiber is of type
$(1,1)$. Demailly--P\u{a}un's deformation theorem for the K\"ahler cone
\cite[Theorem~0.9]{DP04} therefore identifies the K\"ahler cones in the
flat marking outside a countable union $E$ of proper analytic subsets of
$\Delta^*$. Here the projected Gauss--Manin connection is the full flat
connection because $H^{2,0}=H^{0,2}=0$. The complement of $E$ is path
connected: two points can be joined by a two-segment polygonal path with
an intermediate point outside the countably many real lines that would
force either segment to meet $E\cup\{0\}$.

Choose a rational K\"ahler class at some $t_0\notin E$ and clear its
denominators. It gives a fixed integral class on $X$, whose restrictions
are K\"ahler at every $t\notin E$. Proposition~\ref{prop:extendclass}
extends it to a line bundle $L$ on $X$. Kodaira's criterion makes $L_t$
ample for $t\notin E$. For every $m\ge1$, Kodaira vanishing and
$K_{X_t}\cong\cO_{X_t}$ give
\[
 h^0(X_t,L_t^{\otimes m})=\chi(X_t,L_t^{\otimes m})
 \qquad(t\notin E).
\]
Riemann--Roch expresses the right-hand side as a polynomial $P(m)$,
independent of $t$, with leading coefficient $c_1(L_t)^n/n!>0$.
Upper semicontinuity, applied for each fixed $m$, gives
$h^0(X_0,L_0^{\otimes m})\ge P(m)$. Thus $L_0$ is big, and $X_0$ is
Moishezon.
\end{proof}

\section{The Beauville--Bogomolov product over the punctured disc}
\label{sec:product-punct}

Throughout this section and the next, we use the following terminology. A
\emph{Beauville--Bogomolov product} is a product of at most one complex torus
and finitely many simply connected factors of two kinds, as in
Theorem~\ref{thm:BB}. An \emph{irreducible holomorphic symplectic factor} is a
compact K\"ahler manifold whose space of holomorphic two-forms is generated by
a symplectic form; a \emph{strict Calabi--Yau factor} is a compact K\"ahler
manifold of dimension at least three with trivial canonical bundle and no
holomorphic forms of intermediate degree. In both cases simple connectedness is
part of the definition.

\subsection{Local product structure}

\begin{lemma}\label{lem:localproduct}
Let $f\colon W\to U$ be a family of compact complex manifolds over a disc
$U$. If $W_u$ is a Beauville--Bogomolov product $S_0\times S_1\times\cdots
\times S_r$, then, after shrinking $U$ about $u$, the family $W$ is isomorphic
over $U$ to the fiber product of deformations of $S_0,S_1,\dots,S_r$.
\end{lemma}

\begin{proof}
Each factor $S_i$ is unobstructed: its Kuranishi family is smooth of dimension
$h^1(S_i,T_{S_i})$. For a torus this is classical, for a strict Calabi--Yau
manifold it is the Bogomolov--Tian--Todorov theorem, and for an irreducible
holomorphic symplectic manifold it is again the Bogomolov--Tian--Todorov
theorem \cite{Tian87,Todorov89}, the $H^0$-vanishing needed for the locally universal property being
$H^0(S_i,T_{S_i})=0$ from \eqref{eq:contraction}.

The tangent bundle of $W_u$ is the direct sum of the pullbacks of the tangent
bundles of the factors. Put $\widehat S_i:=\prod_{j\neq i}S_j$. The K\"unneth formula gives
\begin{equation}\label{eq:kunnethT}
H^1\bigl(W_u,T_{W_u}\bigr)=\bigoplus_{i=0}^{r}H^1(S_i,T_{S_i}).
\end{equation}
Indeed, the summand $H^1(S_i,T_{S_i})\otimes H^0(\widehat S_i,\cO)$ equals
$H^1(S_i,T_{S_i})$ because $h^0(\widehat S_i,\cO)=1$, while the summand
$H^0(S_i,T_{S_i})\otimes H^1(\widehat S_i,\cO)$ vanishes: either
$H^0(S_i,T_{S_i})=0$ for a simply connected factor, or $S_i$ is the torus and
then $\widehat S_i$ is a product of simply connected factors, so that
$H^1(\widehat S_i,\cO)=0$.

The product $W_u$ itself is compact K\"ahler with trivial canonical bundle,
so its deformations are unobstructed as well \cite{Tian87,Todorov89}. Hence the natural map from the
product of the factor deformation bases to the Kuranishi base of $W_u$ is an
isomorphism of germs: both bases are smooth at
the origin, and the tangent map is the isomorphism \eqref{eq:kunnethT}.
Consequently, after shrinking $U$, the family $W$ splits as the fiber product
of the induced factor families over $U$.

We claim that each factor family is again of the same type. Kodaira--Spencer
stability \cite{KS60} shows that all fibers are K\"ahler, and a differentiable
trivialization identifies both the fundamental groups and the restrictions of
$c_1$ of the relative tangent bundle; hence $c_1=0$ for every fiber.

\emph{The torus factor.} Let $q$ be the fiber dimension. A deformation of a
$q$-dimensional complex torus has $b_1=2q$ and carries $q$ independent
holomorphic one-forms. By Yau's theorem \cite{Yau78} it carries a Ricci-flat
K\"ahler metric, with respect to which these forms are parallel by Bochner's
principle. Their wedge is then a nowhere vanishing holomorphic $q$-form
obtained from parallel forms, so the cotangent bundle is trivialized by
parallel forms, the metric is flat, and the fiber is a complex torus.

\emph{The simply connected factors.} For these factors $H^1(\cO)=0$, and $H^2(-,\Z)$ is torsion free by
the universal coefficient theorem. Thus vanishing of the real first Chern
class implies vanishing of the integral first Chern class, and the exponential
sequence gives triviality of the canonical bundle.
For a strict Calabi--Yau factor constancy of the Hodge numbers preserves the
required vanishing of the intermediate holomorphic forms, and the dimension is
constant. For a symplectic factor the relative holomorphic two-forms form a
line bundle on the base by Grauert's theorem, since $h^{2,0}=1$ on every fiber.
A local relative generator is nondegenerate on the central fiber of this
local deformation. Nondegeneracy is open, and properness lets us shrink $U$ so
that it holds everywhere on the family; equivalently, the top exterior power
of the generator is nowhere zero.
\end{proof}

\begin{lemma}\label{lem:factorwise}
A biholomorphism between Beauville--Bogomolov products preserves the torus
factor and acts factor by factor up to a permutation of the simply connected
factors. For an isomorphism between two local product families, the permutation
is locally constant and the induced factor maps are holomorphic over the base.
\end{lemma}

\begin{proof}
We first show that every Ricci-flat K\"ahler metric on one of the simply
connected factors is de Rham irreducible. Its de Rham factors are complete and
simply connected, and there is no Euclidean factor: a flat de Rham factor of a
compact manifold would carry a nonzero global parallel one-form, which is
impossible on a compact manifold of $b_1=0$, since such a form is closed and
would represent a nonzero class in $H^1$. The parallel complex structure
preserves each nonflat de Rham factor, because the holonomy acts on the tangent
spaces of the de Rham factors independently. Each such factor is therefore
K\"ahler and has trivial canonical bundle.

A nontrivial splitting of a strict Calabi--Yau manifold would produce an
intermediate holomorphic form: a holomorphic $p$-form of one factor with
$0<p<\dim$ pulls back to a nonzero holomorphic $p$-form on the product,
contradicting the definition of a strict Calabi--Yau manifold. For an
irreducible holomorphic symplectic manifold the K\"unneth formula, together
with the vanishing of $H^1$ and of holomorphic two-forms other than the period
line, splits its symplectic form as a sum of pullbacks from the factors; the
vanishing of $h^{1,0}$ rules out the mixed terms, nondegeneracy forces a
nonzero two-form on each factor, and $h^{2,0}=1$ then leaves only one factor.
So no simply connected factor splits further.

Choose a product of Ricci-flat K\"ahler metrics on the target product and let
$g$ be its pullback to the source under the given biholomorphism; then $g$ is
Ricci-flat K\"ahler. Since $\pi_1$ of the source is the product of the
fundamental groups of its factors, and all simply connected factors have
vanishing $H^1(-,\R)$, the K\"unneth formula shows that there are no mixed
summands in $H^{1,1}$ and that the K\"ahler class of $g$ is the sum of the
pullbacks of its restrictions to the source factors; each restriction is a
K\"ahler class on its factor. Let $\tilde g$ be the product of the Ricci-flat
K\"ahler metrics on the source factors lying in these restriction classes, which
exist by Yau's theorem. Then $\tilde g$ is Ricci-flat K\"ahler with the same
K\"ahler class as $g$, so $g=\tilde g$ by Yau's uniqueness theorem
\cite{Yau78}.

Thus the given biholomorphism is an isometry between the Riemannian products
$g$ and the chosen target metric. Lift it to the universal covers: by the
uniqueness of the de Rham decomposition it identifies the Euclidean factor and
permutes the irreducible nonflat factors \cite[Theorem 2 and the following
lemma]{Beauville83}. Consequently each component of the original
biholomorphism has derivative zero along the other source factors, and since
the factors are connected, that component depends only on its corresponding
factor.

For families the permutation is locally constant, because the nonzero degree
two cohomology summands pulled back from the simply connected factors are
distinct and form local subsystems in each product chart. Every positive dimensional compact K\"ahler factor has nonzero second
cohomology, because the top power of a K\"ahler class has positive integral.
The factor summands are nonzero and distinct by the K\"unneth decomposition.
Their images under a family isomorphism vary continuously in a finite set of
possible summands, so the permutation is locally constant. An
isomorphism of families induces an isomorphism of the cohomology local systems.
Once the factors are matched, local holomorphic sections of the other factor
submersions give the induced factor maps by restriction, and independence of
the other coordinates makes these maps agree wherever the sections overlap.
\end{proof}

\begin{lemma}\label{lem:finiteorder}
Let $W\to U$ be a family of compact complex manifolds with connected total
space. A finite order holomorphic automorphism over $U$ which fixes one fiber
pointwise is the identity.
\end{lemma}

\begin{proof}
Let $g$ be the automorphism and suppose $g^m=1$. At any point of the fixed
fiber its differential is the identity on the vertical tangent space and on
the quotient coming from $U$. Its matrix is therefore unipotent with a single
possible off-diagonal block. Since $(dg)^m=1$ in characteristic zero, that
block vanishes, and $dg$ is the identity. Restrict the domain of a local coordinate map $z$ to the
intersection of its $g$-translates and average: the map
$\frac1m\sum_{k=0}^{m-1}z\circ g^k$ is $g$-invariant and has invertible
derivative at the fixed point, hence is a coordinate map in which $g$ acts as
the identity. It follows that $g$ is the identity on a neighbourhood of that
point and, $W$ being connected, on all of $W$.
\end{proof}

\subsection{The product cover}

\begin{theorem}\label{thm:productcover}
Let $\pi\colon X\to\Delta$ be a family of compact complex manifolds. Suppose
that $X_t$ is K\"ahler and $c_1(X_t)=0$ in $H^2(X_t,\R)$ for every
$t\in\Delta^*$. Then there are a finite \'etale cover $p\colon\tX\to X$ such
that $\tpi=\pi\circ p$ has connected fibers, a finite set $I=\{i_1,\dots,i_s\}$,
and families $Q^*_i\to\Delta^*$ with connected fibers, together with an
isomorphism
\[
\Phi^*\colon\tX|_{\Delta^*}\xrightarrow{\ \sim\ }
Q^*_{i_1}\times_{\Delta^*}\cdots\times_{\Delta^*}Q^*_{i_s}
\]
over $\Delta^*$. The torus factor, if present, has index $0$, and each other
factor has fibers of one of the two simply connected types.
\end{theorem}

\begin{proof}
Fix $t^*\in\Delta^*$ and choose a connected finite \'etale
Beauville--Bogomolov cover of $X_{t^*}$; let $H\subseteq\pi_1(X_{t^*})$ be the
defining subgroup. The group $H$ is free abelian of finite rank: it is the fundamental group of
the chosen product cover, whose only nontrivial fundamental group comes from
the torus factor.
Ehresmann's theorem \cite{Ehresmann47} and the contractibility of $\Delta$ show
that the inclusions of the fibers induce isomorphisms on fundamental groups and
on the corresponding local systems. The subgroup $H$ therefore defines a finite
topological cover $\tX\to X$ whose restrictions to the fibers are connected; the
lifted complex structure makes it a finite \'etale cover over the entire disc.

Let $V\subset\Delta^*$ be the locus where the covered fiber $\tX_t$ is itself a
Beauville--Bogomolov product. It contains $t^*$, and it is open by
Lemma~\ref{lem:localproduct}. To prove that $V$ is closed in $\Delta^*$, fix
$b\in\overline V\cap\Delta^*$ and choose a Beauville--Bogomolov cover of
$\tX_b$. This cover extends over a small disc $U$ around $b$ to a connected
finite \'etale cover $W\to\tX|_U$ with connected fibers, by the same topological
argument. The cover is Galois: its defining subgroup $N\subseteq H$ is normal
because $H$ is abelian. Denote its deck group by $G=H/N$.

After shrinking $U$, Lemma~\ref{lem:localproduct} identifies $W$ with
$A\times_US$, where $A$ is the torus family --- possibly a family of points ---
and $S$ is the product of the simply connected factor families.
Lemma~\ref{lem:factorwise} gives separate $G$-actions on $A$ and $S$. The
action on each $\pi_1(W_u)\cong N$ is induced by conjugation in $H$ and hence is
trivial. Since the inclusion of the torus factor induces an isomorphism on
fundamental groups, every induced automorphism of $A_u$ is a translation.

Now choose $u\in U\cap V$; this set is nonempty because $U$ meets $V$ and $V$
is open. Every connected finite cover of $\tX_u$ comes from a cover of its
torus factor, by the classification of covering spaces and the simple
connectedness of the other factors. In that splitting the deck group acts
trivially on the simply connected part. By Lemma~\ref{lem:factorwise} the same
holds in the splitting of $W$ inherited from $W\to\tX|_U$.
Lemma~\ref{lem:finiteorder} applied to $S\to U$ then makes the action on $S$
trivial throughout $U$. The quotient fiber at $b$ is therefore the product of
$A_b/G$, a torus, and $S_b$. This proves that $V$ is closed in the connected
punctured disc, so $V=\Delta^*$.

It remains to glue the local products. By Lemma~\ref{lem:factorwise} the only
possible monodromy of the factor labels is a permutation of the simply connected
factors, and such a permutation permutes the distinct nonzero summands of second
cohomology pulled back from these factors. The local system
$R^2\tpi_*\R$ is constant on the whole disc, so this permutation is trivial.
After choosing global labels, the factor maps on overlaps satisfy the cocycle
condition because their product is the transition isomorphism. They glue the
factor families and their product identifications over $\Delta^*$. The glued
families are proper holomorphic submersions since these properties are local on
the base.
\end{proof}

\section{Extension of the relative product}\label{sec:product}

\subsection{Triviality of the relative canonical bundle}

\begin{lemma}\label{lem:trivialcanonical}
Let $f\colon\tX\to\Delta$ be a family of compact complex manifolds and assume
that $K_{\tX_t}\cong\cO_{\tX_t}$ for $t\neq0$. Then $K_{\tX/\Delta}$ is
holomorphically trivial.
\end{lemma}

\begin{proof}
The direct image $f_*K_{\tX/\Delta}$ is coherent by Grauert's theorem
\cite{Grauert60} and torsion free, since a section of a
line bundle which vanishes on a nonempty open set vanishes identically. It has
rank one: by Grauert's theorem and base change over $\Delta^*$, its fibers
there are $H^0(\tX_t,K_{\tX_t})\cong\C$ for $t\neq0$. Since $\Delta$ is a
smooth curve and $f_*K_{\tX/\Delta}$ is coherent and torsion free, it is
locally free of rank one, hence a line bundle; it is trivial because
$\Pic(\Delta)=0$. A generator gives a relative top form $\tau$ which is nowhere
zero on every punctured fiber, so that $\tau$ trivializes $K_{\tX/\Delta}$ over
$\Delta^*$.

The divisor of $\tau$ is supported on $\tX_0$. Since that fiber is smooth
and connected, it is irreducible, and the divisor is $m\tX_0$ for an integer
$m\geq0$. We show that $m=0$. Suppose otherwise.
Since $f_*K_{\tX/\Delta}\cong\cO_\Delta$ is generated by $\tau$, every global
relative top form is $\varphi\tau$ for a holomorphic function $\varphi$ on
$\Delta$. Now $\tau$ vanishes to some finite order $m\ge1$ along $\tX_0$, so
$\tau/t^m$ is a relative top form defined near $\tX_0$ and holomorphic there; by
the removable singularity theorem this extends to a global relative top form,
hence equals $\psi\tau$ for some $\psi$ holomorphic on $\Delta$. Then
$1=t^m\psi$ on $\Delta^*$, impossible since the left hand side does not vanish
at $0$. Thus $m=0$. A nonzero holomorphic section of a line bundle on a complex
manifold has a zero locus of pure codimension one if it has any zeros. Hence
$\tau$ has no zeros at all and trivializes $K_{\tX/\Delta}$.
\end{proof}

\subsection{Saturated periods}

Fix the cover and the product $\Phi^*$ of Theorem~\ref{thm:productcover}. The
argument of that theorem applies verbatim to any fixed finite \'etale cover with
connected fibers and a prescribed Beauville--Bogomolov product over the
puncture; we use only that the covered punctured fibers are K\"ahler with
holomorphically trivial canonical bundle. In particular, when the original
punctured fibers only have vanishing real first Chern class,
Theorem~\ref{thm:BB} shows that the product cover has the required trivial
canonical bundles.

Let $I$ be the set of positive dimensional factor labels, with $0$ reserved for
the torus when present. Write $f\colon\tX\to\Delta$ for the family on the cover
and $V^k_i$ for the flat K\"unneth summand pulled back from degree $k$ of the
factor $i$. Let $d_i$ be the dimension of that factor and $n=\sum_{i\in I}d_i$
the relative dimension. We omit factors of dimension zero throughout. The
labels and the K\"unneth algebra are marked by the differentiable
trivialization over the disc. In particular, the cup product of nonzero classes
in the top degree summands $V^{d_i}_i$ is nonzero. In this section all markings
and forms are on the cover; we write $\int_i$ for evaluation on the marked
fundamental class of the factor $i$.

\begin{lemma}\label{lem:saturated}
For every factor the marked top form line extends holomorphically over zero.
For any nonzero generator $u_i(0)$ of the limiting line, the factor pairing
$\int_i u_i(0)\wedge\overline{u_i(0)}$ is nonzero. The degree two Hodge
splitting of the punctured fibers extends as a smooth splitting of the
complexified flat cohomology bundle over the disc.
\end{lemma}

\begin{proof}
\emph{Step 1: the top form line.} For each $i$ let $L_i$ be the kernel of
\[
f_*\Omega^{d_i}_{\tX/\Delta}\longrightarrow
\bigl(H^{d_i}_\C/V^{d_i}_i\bigr)\otimes\cO_\Delta ,
\]
the map sending a relative form to its fiberwise de Rham class modulo
$V^{d_i}_i$. Relative sections are fiberwise closed by holomorphic continuation
from the K\"ahler puncture, and the relative holomorphic de Rham comparison
gives a holomorphic period map. The direct image is torsion free on the smooth
curve $\Delta$, hence locally free. The kernel is saturated because the quotient
embeds in a vector bundle. It has rank one by the Hodge--K\"unneth decomposition
on $\Delta^*$ and the trivial canonical bundle of each factor.

Choose a generator $\theta_i$ of $L_i$ near zero. Its restriction to $\tX_0$ is
a nonzero form: otherwise $\theta_i$ would be divisible by $t$ in the direct
image, and saturation would make it divisible by $t$ in $L_i$. Its period,
however, may a priori vanish at zero. Write
\[
[\theta_{i,t}]=t^{k_i}u_i(t),\qquad k_i\ge0,\quad u_i(0)\neq0,\quad
u_i(t)\in V^{d_i}_i .
\]
The line generated by $u_i$ is the saturation of the period image in the flat
bundle.

By Lemma~\ref{lem:trivialcanonical} there is a nowhere vanishing relative top
form $\tau$ on the cover. The product of the $\theta_i$ satisfies
\[
\bigwedge_{i\in I}\theta_i=t^{\ell}b(t)\,\tau,\qquad b(0)\neq0 ,
\]
where the coefficient is a function of the base by connectedness and
compactness of the fibers. The class $[\tau_0]$ is nonzero, since
$i^{n^2}\int_{\tX_0}\tau_0\wedge\overline{\tau_0}>0$ by the volume
form property and Stokes' theorem.
The cup product $\prod_iu_i(0)$ is nonzero by the marked K\"unneth formula.
Comparison of vanishing orders in marked cohomology gives $\ell=\sum_ik_i$ and
$\prod_iu_i(0)=b(0)[\tau_0]$. Consequently
\[
\prod_{i\in I}\int_i u_i(0)\wedge\overline{u_i(0)}
=\varepsilon\,|b(0)|^2\int_{\tX_0}[\tau_0]\wedge\overline{[\tau_0]} ,
\]
where $\varepsilon=\pm1$ is the K\"unneth sign. The right hand side is nonzero
because $\tau_0$ is a volume form, so every factor pairing is nonzero,
independently of the orders $k_i$.

\emph{Step 2: the torus factor.} For a torus factor, saturate the rank $d_i$
period image in $V^1_i\otimes\cO_\Delta$; equivalently, take the saturation of
the image of the kernel of the period map to the quotient by $V^1_i$.
Saturation of a subsheaf of a locally free sheaf on a smooth curve is a
subbundle. Its determinant in $\bigwedge^{d_i}V^1_i$ is the top form line on
the puncture, hence also at zero. Let $W\subset V^1_i$ be the value of this rank
$d_i$ subbundle at $0$, and write $\overline W$ for the value of its conjugate.
A generator of $\bigwedge^{d_i}W$ differs from the limiting top form
$u_i(0)$ by a nonzero scalar. Hence the nonvanishing of the
factor pairing $\int_iu_i(0)\wedge\overline{u_i(0)}$ says exactly that
$\bigwedge^{d_i}W\wedge\bigwedge^{d_i}\overline W\neq0$ in the one dimensional
space $\bigwedge^{2d_i}V^1_i$. This in turn holds if and only if
$W\cap\overline W=0$: a nonzero intersection would allow a common vector to be
completed to bases of both $W$ and $\overline W$, making the product of the two
determinant lines vanish. Therefore $V^1_i=W\oplus\overline W$, and taking
exterior powers extends the Hodge splitting to $V^2_i$.

\emph{Step 3: the symplectic factor.} For an irreducible holomorphic symplectic
factor of dimension $2m_i$, saturate the period line in $V^2_i$, and write
$v_i(t)$ for its nonvanishing holomorphic frame and $q_i$ for the marked
Beauville--Bogomolov form. Use the Fujiki relation
$\int_ix^{2m_i}=c_iq_i(x)^{m_i}$ in the marked cohomology algebra
\cite{Fujiki87}. Since $q_i(v_i(t))=0$ on the puncture, $q_i(v_i(0))=0$ by
continuity. Applying the Fujiki relation to the class $v_i(0)+sw$ and
polarizing, the coefficient of $s^{m_i}$ equals
\[
\binom{2m_i}{m_i}\int_iv_i(0)^{m_i}w^{m_i}
=c_i\bigl(2q_i(v_i(0),w)\bigr)^{m_i}
\]
for every $w\in H^2$. Suppose that $v_i(0)^{m_i}=0$ as a class. Then the left
hand side vanishes for every $w$, hence $q_i(v_i(0),w)=0$ for every $w$, and the
nondegeneracy of $q_i$ forces $v_i(0)=0$, contradicting the fact that
$v_i(0)$ is the value at zero of a nonvanishing frame. Therefore
$v_i(0)^{m_i}\neq0$, and $v_i(0)^{m_i}$ spans the limiting top form line.

Next, polarize the Fujiki relation in two classes $x,y$ with $q_i(x)=q_i(y)=0$:
applying it to $x+sy$ and comparing coefficients of $s^{m_i}$ gives
\[
\int_ix^{m_i}y^{m_i}
=\frac{c_i\,2^{m_i}}{\binom{2m_i}{m_i}}\,q_i(x,y)^{m_i}.
\]
Taking $x=v_i(0)$ and $y=\overline{v_i(0)}$, whose $q_i$-square vanishes because
$q_i$ is defined over $\R$ and $q_i(v_i(0))=0$, and using the nonvanishing of
the factor pairing established above, we obtain
$q_i(v_i(0),\overline{v_i(0)})\neq0$. It is positive by continuity, since
$q_i(v_i(t),\overline{v_i(t)})>0$ on the puncture. Consequently the $(2,0)$
line, the $(0,2)$ line and their $q_i$-orthogonal complement are complementary
at zero --- the restriction of $q_i$ to the span of the first two lines has Gram
matrix $\bigl(\begin{smallmatrix}0&c\\ c&0\end{smallmatrix}\bigr)$ with
$c\neq0$, hence is nondegenerate --- and the last of these spaces extends the
$(1,1)$ summand.

\emph{Step 4: the full cover.} The strict Calabi--Yau factors have all of $H^2$
of type $(1,1)$, because $h^{2,0}=h^{0,2}=0$ there. Since the nontorus factors have vanishing $H^1$, there are no mixed degree
two K\"unneth terms. This proves
the assertion for the full cover. The extended splitting at zero is defined
inside the marked cohomology bundle, and its construction uses the Hodge
decomposition only on $\Delta^*$.
\end{proof}

\begin{lemma}\label{lem:movingsection}
There is a smooth section $a(t)$ of the marked real degree two cohomology of the
cover, of type $(1,1)$ on the puncture and K\"ahler on a dense open subset of
it, with constant positive top self-intersection.
\end{lemma}

\begin{proof}
Use Lemma~\ref{lem:saturated} to extend the projected Gauss--Manin connection
$\nabla^{1,1}$ on the real limiting $(1,1)$ bundle smoothly over zero: the
splitting supplies a smooth projection onto the $(1,1)$ part, and the connection
is then the composition of the Gauss--Manin connection with that projection.

Apply Demailly--P\u{a}un's theorem on the transport of the K\"ahler cone
\cite[Theorem 0.9]{DP04} to the family restricted to $\Delta^*$, where every
fiber is K\"ahler. Outside a countable set $E\subset\Delta^*$, parallel
transport along paths in $\Delta^*\setminus E$ preserves the K\"ahler cones.

Choose $t^*\notin E\cup\{0\}$ and a K\"ahler class $a^*$ there. Define $a(t)$ by
parallel transport for the extended connection along the straight segment
from $t^*$ to $t$. Smooth dependence of
parallel transport on the endpoint gives a smooth section on the entire disc,
including zero. For all $t$ outside the union of the real lines through $t^*$
and the points of $E\cup\{0\}$, the segment stays in the good locus, so $a(t)$
is K\"ahler. There are countably many such lines, so the remaining endpoints are
dense by the Baire category theorem.

Along each of these segments the ordinary Gauss--Manin derivative of $a$ has only
$(2,0)$ and $(0,2)$ components, because $a(t)$ is of type $(1,1)$ along the
segment. Hence the derivative of the top self-intersection $\int a(t)^n$ is zero
by type reasons, and continuity makes the top self-intersection constant on the
disc, equal to the positive volume at $t^*$.

The K\"ahler locus is open on the puncture. Near a K\"ahler parameter choose
smooth auxiliary K\"ahler metrics by K\"ahler stability \cite{KS60}. Their
fiberwise harmonic representatives of $a(t)$ are smooth and, being of type
$(1,1)$, are the $(1,1)$-parts of these representatives; at the given parameter
the $\pdbar$-lemma expresses the difference from a positive representative by a
potential. Extend that potential smoothly to neighboring fibers and use openness
of positivity.
\end{proof}

\subsection{Removal of poles of the factor forms}

Normalize the factor forms by $\Omega_i=t^{-k_i}\theta_i$, so that their periods
are the holomorphic sections $u_i$ of Lemma~\ref{lem:saturated}. For the
auxiliary periods we use the same kernel of the period map to the quotient by
the corresponding K\"unneth summand. In degree one for the torus and in degree
two for a symplectic factor this kernel maps isomorphically onto its image over
$\Delta^*$ by K\"ahler Hodge theory; its image and its saturation agree after
inverting $t$. A holomorphic frame of the saturation therefore has meromorphic
relative lifts with possible poles only at zero. All the chosen forms are
fiberwise closed on the puncture and are pulled back from the designated
factors.

Let $U_a$ be the K\"ahler locus of the moving class of
Lemma~\ref{lem:movingsection}. For $t\in U_a$ let $\omega_t$ be the Ricci-flat
representative of $a(t)$ given by Yau's theorem \cite{Yau78}, and let $h_t$ be a
smooth family of Gauduchon forms through zero, which exists by
Lemma~\ref{lem:parametricGauduchon}. Proposition~\ref{prop:mass} applies to the
smooth marked section $a(t)$ and, since $\omega_t\wedge h_t^{n-1}$ is a positive
$(n,n)$-form, gives
\begin{equation}\label{eq:masstr}
\int_{\tX_t}\tr_{h_t}\omega_t\,h_t^{\,n}
=n\int_{\tX_t}\omega_t\wedge h_t^{n-1}\le C\qquad(t\in U_a),
\end{equation}
where we used $(\tr_{h_t}\omega_t)h_t^n=n\,\omega_t\wedge h_t^{n-1}$. The bound
is uniform as $t\to0$ in $U_a$.

Let $\eta$ be any of the meromorphic relative forms chosen above, of degree
$p\ge1$. On $U_a$ the Bochner principle makes $\eta_t$ parallel with respect to
the Ricci-flat metric $\omega_t$. Hence its pointwise squared norm for
$\omega_t$ is a dimensional constant times
\[
\frac{\bigl|i^{p^2}\int_{\tX_t}[\eta_t]\wedge\overline{[\eta_t]}\wedge
a(t)^{n-p}\bigr|}{\int_{\tX_t}a(t)^n}.
\]
The numerator is bounded in the marked cohomology algebra and the denominator is
a fixed positive number. Comparison of norms in a frame diagonalizing $\omega_t$
relative to $h_t$ gives
\[
|\eta_t|_{h_t}^{2/p}\le C_{n,p}\,|\eta_t|_{\omega_t}^{2/p}\,
\tr_{h_t}\omega_t .
\]
Together with \eqref{eq:masstr} this bounds the integral of $|\eta_t|_{h_t}^{2/p}$
uniformly on $U_a$ near zero.

If $\eta=t^{-m}\xi$ had a genuine pole of order $m>0$, its holomorphic leading
coefficient $\xi|_{\tX_0}$ would be nonzero on a smaller relative coordinate
box. On that box $|\xi_t|_{h_t}$ and the $h_t$-volume of the fiber slice are
bounded below, so the integral would grow at least like a positive constant
times $|t|^{-2m/p}$ as $t\to0$ in $U_a$. This contradiction shows that every
chosen form, including the auxiliary one-forms and two-forms, extends
holomorphically.

For the factor volume forms this gives $k_i=0$: indeed $\theta_i|_{\tX_0}\neq0$
as a form while $[\theta_{i,t}]=t^{k_i}u_i(t)$ with $u_i(0)\neq0$, so $k_i>0$
would exhibit a pole of the extension, which we just excluded. In particular the
central periods are $u_i(0)\neq0$. The comparison of vanishing orders in
Lemma~\ref{lem:saturated} then gives
\begin{equation}\label{eq:wedgeproduct}
\bigwedge_{i\in I}\Omega_i=b(t)\,\tau,\qquad b(0)\neq0 .
\end{equation}
After shrinking, the wedge \eqref{eq:wedgeproduct} is nowhere zero. The extended
forms remain fiberwise closed, so their kernels are integrable, and the
nonvanishing of the wedge gives the rank statements used in the proof of
Theorem~\ref{thm:product} below. The Pl\"ucker relations for decomposability hold on the puncture and extend
by holomorphic continuation. Together with \eqref{eq:wedgeproduct}, they imply
that each extended $\Omega_i$ is a decomposable form of constant rank $d_i$.

On the torus, the wedge of the extended one-forms equals a holomorphic unit
times the extended factor volume form. On an irreducible symplectic factor of
dimension $2m_i$, the $m_i$-th power of the extended two-form satisfies the same
assertion. Indeed, the corresponding identity of period sections holds with a
unit coefficient, by the determinant description for the torus and the Fujiki
identity for the symplectic factor. Injectivity of holomorphic forms into de
Rham cohomology on the punctured fibers gives the identity of forms there, and
holomorphic continuation gives it on the central fiber.

\subsection{The product across the central parameter}

\begin{theorem}\label{thm:product}
Let $\tX\to\Delta$ be a family of compact complex manifolds with a prescribed
Beauville--Bogomolov product $\Phi^*$ over $\Delta^*$ as above, with factor set
$I$. Then there are families $Q_i\to\Delta$ with connected fibers and an
isomorphism
\[
\Phi\colon\tX\xrightarrow{\ \sim\ }\prod_{\Delta,\,i\in I}Q_i
\]
over $\Delta$. For each $i$ there is an isomorphism $Q_i|_{\Delta^*}\cong
Q^*_i$ which identifies the restriction of the $i$-th projection of $\Phi$ with
the $i$-th projection of $\Phi^*$. The construction takes place on the given
cover and does not require a hypothesis on the central fiber.
\end{theorem}

\begin{proof}
For an empty product the family is the base, and for a single factor we take
$Q_i=\tX$. For more than one factor, work first over a small disc around zero on
which the relative forms of the preceding subsection are defined and
\eqref{eq:wedgeproduct} is nowhere zero.

In the relative tangent bundle put
\[
F_i=\bigcap_{j\neq i}\ker\Omega_j,\qquad
G_i=\ker\Omega_i=\bigoplus_{j\neq i}F_j,
\]
where $\ker\Omega_i$ denotes the kernel of contraction by $\Omega_i$. The
decomposability of the $\Omega_i$ extends from the puncture, and the nonvanishing
wedge \eqref{eq:wedgeproduct} gives constant ranks and the displayed direct sum.
The forms are closed along the fibers, so their kernels and the displayed
intersections are integrable. The foliations are vertical, and on the puncture
their leaves are the indicated products of factors with the other coordinates
held fixed.

Fix a smooth Hermitian background $h$ and work over a compact subdisc. Set
$\Xi_i=\bigwedge_{j\in I\setminus\{i\}}\Omega_j$ and $e_i=n-d_i$. On the
distribution $G_i$ the fiberwise closed form $i^{e_i^2}\Xi_i\wedge
\overline{\Xi_i}$ is a positive leafwise volume form: restricted to a leaf of
$G_i$ it is a nonzero decomposable form of top degree $2e_i$, because on the
puncture the leaves are products of the factors $j\neq i$ and $\Xi_i$ restricts
to the product of their top forms. By compactness of the subdisc it is uniformly
comparable to the leafwise volume density of $h$, in the sense that the two are
bounded above and below by positive constants times each other.

We now bound the leaf volumes. Over the puncture the family is a product, so
a $G_i$-leaf in $\tX_t$ is a product of the factors with labels $j\neq i$.
Its homology class is independent of the fixed point in the remaining factor
and is locally constant under the marking. Since the local system on the disc
is constant, this is one fixed marked class. Pairing it with the smooth
family of de Rham classes $[i^{e_i^2}\Xi_i\wedge\overline{\Xi_i}]$ --- which is
smooth across zero by the extension of the forms --- bounds the $h$-volumes of
all punctured $G_i$-leaves by a constant independent of $t$.

Let $x\in\tX_0$ and choose $x_\nu\in\tX_{t_\nu}$ with $x_\nu\to x$ and
$t_\nu\neq0$. The compact $G_i$-leaves through $x_\nu$ have bounded volume and
supports in a fixed compact subset of the total space. Bishop's compactness
theorem \cite{Bishop64}, in the cycle-theoretic form of
\cite[Theorem~1]{Barlet99}, gives a subsequence converging to a compact cycle
supported on $\tX_0$ and containing $x$. In a relative foliation box write the
coordinates as $(z,w,t)$ with $z$ along the plaques. Each approximating leaf
support contains the whole connected plaque through each of its points, so local
convergence of supports gives the same saturation for the limit support.
Intersect this support with a transverse slice $z=z_0$ in the central fiber. The
intersection is analytic and zero dimensional, since a positive dimensional
component would give a component of the support of dimension exceeding the rank
of $G_i$. The support is thus locally a disjoint union of plaques, and its
connected component through $x$ is exactly the central leaf through $x$, which
is therefore compact.

Every $G_i$-leaf has trivial holonomy in the total space. For a loop on a
central leaf, choose its local transversal by taking a small $F_i$-plaque
through a local section over each parameter; this is a holomorphic transversal
because $F_i$ and $G_i$ are complementary in the vertical tangent space. The
holonomy germ preserves the parameter and relates points on the same $G_i$-leaf.
For each punctured parameter the transversal lies in a single $F_i$-leaf, and in
the punctured product a $G_i$-leaf and an $F_i$-leaf meet in exactly one point.
The germ is therefore the identity on the puncture, and holomorphic continuation
gives the identity at zero. The same product description gives trivial holonomy
for the punctured leaves.

Reeb stability for a compact leaf with trivial holonomy gives a saturated
neighborhood whose leaves each meet a chosen local transversal exactly once
\cite[Theorem 2.9]{Moerdijk03}. The transversal can be chosen holomorphic,
because the foliation is holomorphic. The resulting quotient chart is then
holomorphic: changes of transversal are holomorphic holonomy maps. These
saturated neighborhoods can be chosen inside any prescribed neighborhood of
the compact leaf; no holomorphic triviality of the family of leaves is asserted. Distinct leaves are disjoint compact subsets of the Hausdorff
manifold $\tX$, so they have disjoint neighborhoods and therefore disjoint
saturated neighborhoods after this shrinking. The leaf space
$Q_i=\tX/G_i$ is consequently Hausdorff. It is
second countable because the quotient map is open. The local holomorphic
transversals and holomorphic transition maps make it a complex manifold, and the
quotient map is a holomorphic submersion. The projection $Q_i\to\Delta$ is a
submersion because the original base projection is one and is constant along the
leaves. Its fibers are connected, being continuous images of the connected
fibers of $f$. It is proper because the inverse image of a compact subset of the disc
is the image of the corresponding compact subset of $\tX$.

Write $q_i$ for the quotient maps and $p^*_i$ for the projections of the
prescribed $\Phi^*$. Over the puncture of the small disc, $q_i$ and $p^*_i$ have
the same connected fibers. Local holomorphic sections of both submersions give a
unique holomorphic isomorphism $\psi_i$ of their targets with $\psi_i\circ
q_i=p^*_i$ there. The product of the quotient maps has invertible differential
on each fiber, by the complementary tangent decomposition, and preserves the base
coordinate. It is proper: if $K$ is compact in the target and $B$ is its compact image
in the base, then its inverse image is a closed subset of the compact space
$f^{-1}(B)$. A proper local biholomorphism has open and closed image, hence is
a covering of the connected target. Its degree is one on the puncture and
therefore everywhere. This gives the biholomorphism near zero. Gluing the factors to $Q^*_i$
via the maps $\psi_i$ gives the asserted families and product isomorphism over
the entire disc.
\end{proof}

\section{Irreducible holomorphic symplectic limits}\label{sec:ihs}

This section is devoted to the following statement.

\begin{theorem}\label{thm:ihs}
Let $\pi\colon X\to\Delta$ be a proper holomorphic submersion with connected
fibers. Suppose that $X_t$ is an irreducible compact hyperk\"ahler manifold for every
$t\in\Delta^*$. Then $X_0$ is bimeromorphic to an irreducible compact
hyperk\"ahler manifold. In particular $X_0$ is simply connected, admits a
nondegenerate holomorphic two-form, and satisfies $h^{2,0}(X_0)=1$.
\end{theorem}

\subsection{Approximation with Picard number zero}

We first fix the distances used in the approximation argument. Fix a
differentiable marking and put $V=H^2(M,\R)$, where $M$ is the underlying
smooth manifold of a fiber. Choose once and for all a real basis
$e_1,\dots,e_b$ of $V$ and set
\[
\Bigl\|\sum_{\nu=1}^{b}x_\nu e_\nu\Bigr\|_V=\Bigl(\sum_{\nu=1}^{b}
x_\nu^2\Bigr)^{1/2},\qquad d_V(c,c')=\|c-c'\|_V .
\]
Fix also a smooth Riemannian metric $h$ on $M$, with Levi-Civita connection
$\nabla$, and use the tensor norms
$\|T\|_{C^r(h)}=\max_{0\le j\le r}\sup_M|\nabla^jT|_h$. Hermitian metrics are
regarded as their fundamental real two-forms. All comparisons of classes and
tensors use these fixed choices and the indicated differentiable
trivializations.

If $Z$ is compact K\"ahler, $g$ a Gauduchon form on $Z$, and
$c\in H^{1,1}(Z,\R)$, we write
\[
\mu_{Z,g}(c)=\int_Z\eta\wedge g^{d-1},\qquad d=\dim_\C Z,
\]
where $\eta$ is any smooth closed real $(1,1)$-representative of $c$. The
$\pdbar$-lemma and $\pdbar g^{d-1}=0$ show that this is independent of $\eta$;
the same pairing is obtained from any closed real $(1,1)$-current representing
$c$, by the $\pdbar$-lemma for currents and integration by parts.

\begin{lemma}\label{lem:perturb}
Let $Z=(M,J)$ be an irreducible compact hyperk\"ahler manifold of complex
dimension $d$, let $g$ be a smooth Gauduchon form, and let $a\in\cC^+_Z$
satisfy $q(a)=1$. For every $\varepsilon>0$ and every integer $r\ge0$ there are
a marked small deformation $Z'=(M,J')$, a smooth Gauduchon form $g'$ on $Z'$,
and a K\"ahler class $a'$ such that
\begin{align}
&\rho(Z')=0,\qquad q(a')=1,\qquad d_V(a',a)<\varepsilon,\label{eq:pert1}\\
&\|J'-J\|_{C^r(h)}+\|g'-g\|_{C^r(h)}<\varepsilon,\label{eq:pert2}\\
&|\mu_{Z',g'}(a')-\mu_{Z,g}(a)|<\varepsilon.\label{eq:pert3}
\end{align}
Moreover the deformation $Z'$ may be required to carry fiberwise holomorphic
charts on $M$ which are $C^r$-close to prescribed ones, in the sense made
precise in the proof.
\end{lemma}

\begin{proof}
Take the full local deformation family $\cZ\to B$ of $Z$ with $\cZ_0=Z$. After
shrinking, $B$ is a coordinate ball, the fibers are irreducible compact
hyperk\"ahler manifolds by Theorem~\ref{thm:localTorelli}, and a differentiable
trivialization gives complex structures $J_s$ on $M$ with $J_0=J$. The markings
identify their second cohomology with $V$ and their Beauville--Bogomolov forms
with $q$. Choose a nonvanishing local holomorphic lift $v(s)\in V_\C$ of the
period line and put $Q(s)=q(v(s),\overline{v(s)})>0$.

Define the real class
\begin{equation}\label{eq:bs}
b(s)=a-\frac{q(a,v(s))}{Q(s)}\,\overline{v(s)}
-\frac{q(a,\overline{v(s)})}{Q(s)}\,v(s).
\end{equation}
Since $q(v(s))=q(\overline{v(s)})=0$, a direct calculation gives
$q(b(s),v(s))=q(b(s),\overline{v(s)})=0$. By \eqref{eq:qHodge} the class $b(s)$
therefore has type $(1,1)$ on $Z_s$. Moreover $b(0)=a$, so $q(b(s))>0$ after
shrinking $B$. Consequently
\begin{equation}\label{eq:As}
A(s)=\frac{b(s)}{\sqrt{q(b(s))}}
\end{equation}
is smooth, satisfies $A(0)=a$, and has square one.

Choose a smooth local family of K\"ahler forms $\omega_s$ by K\"ahler stability
\cite{KS60}. At $s=0$ the classes $a$ and $[\omega_0]$ lie in the same positive
cone component, whence $q(a,[\omega_0])>0$; by continuity, shrink $B$ so that
$q(A(s),[\omega_s])>0$. The Lorentzian signature on $H^{1,1}(Z_s,\R)$ then gives
$A(s)\in\cC^+_{Z_s}$ for every $s\in B$.

We next extend $g$ while retaining its exact value at $s=0$: by
Lemma~\ref{lem:parametricGauduchon} there is a smooth family $G_s$ of Gauduchon
forms on the fibers of $\cZ\to B$ with $G_0=g$.

Let $\eta_s$ be the de Rham harmonic representative of $A(s)$ for the K\"ahler
metric $\omega_s$. K\"ahler Hodge theory gives a closed real $(1,1)$-form, and
parameter dependent elliptic theory gives smooth dependence on $s$. It follows
that
\[
F(s):=\int_M\eta_s\wedge G_s^{d-1}=\mu_{Z_s,G_s}\bigl(A(s)\bigr)
\]
is smooth with $F(0)=\mu_{Z,g}(a)$. Therefore the set
\[
U_{\varepsilon,r}=\left\{s\in B:\begin{aligned}
&d_V(A(s),a)<\varepsilon,\\
&\|J_s-J\|_{C^r(h)}+\|G_s-g\|_{C^r(h)}<\varepsilon,\\
&|F(s)-F(0)|<\varepsilon
\end{aligned}\right\}
\]
is an open neighbourhood of $0$.

By Corollary~\ref{cor:density} the set
$B^{\mathrm{gen}}=\{s\in B:\rho(Z_s)=0\}$ is dense in $B$. Choose
$s\in U_{\varepsilon,r}\cap B^{\mathrm{gen}}$ and set
\[
Z'=Z_s,\qquad g'=G_s,\qquad a'=A(s).
\]
Then $a'\in\cC^+_{Z'}=\cK(Z')$ by \eqref{eq:cone}, and the three error bounds
follow from the definition of $U_{\varepsilon,r}$. Note that a path in $B$ from
$0$ to $s$ preserves the connected marked deformation component.

It remains to justify the assertion about charts, which is needed in
Proposition~\ref{prop:approx}. The local coordinate theorem with parameters
supplies smoothly varying fiberwise holomorphic charts on relatively compact
smaller coordinate domains. For a finite collection of such domains,
$C^r$-closeness of their coordinate maps to the base point maps is a finite
collection of additional open conditions on $s$, so the same choice of $s$
applies to all of them.
\end{proof}

\begin{proposition}\label{prop:approx}
Let $\pi\colon X\to\Delta$ be a proper holomorphic submersion with connected
fibers. Suppose that $X_t$ is an irreducible compact hyperk\"ahler manifold for every
$t\in\Delta^*$. Then there exist marked irreducible compact hyperk\"ahler
manifolds $Z_i=(M,J_i)$, smooth Gauduchon forms $g_i$, and K\"ahler classes
$a_i$ such that
\[
J_i\longrightarrow J_0,\quad g_i\longrightarrow g_0\ \text{ in }C^\infty,
\qquad \rho(Z_i)=0,
\]
\[
a_i\longrightarrow a_0,\qquad q(a_i)=q(a_0)=1,
\]
and every closed positive $(1,1)$-current $T_i$ in $a_i$ satisfies
\begin{equation}\label{eq:massbound}
\int_{Z_i}T_i\wedge g_i^{d-1}\le C .
\end{equation}
All marked $Z_i$ lie in the same connected marked deformation component as the
punctured fibers.
\end{proposition}

\begin{proof}
Let $g_t$ be a smooth family of Gauduchon forms on the original disc, which
exists by Lemma~\ref{lem:parametricGauduchon}. Proposition~\ref{prop:mass}
bounds the class pairing $m(t)$ of $a_t$ by a common constant; here $a_t$ is the
moving K\"ahler class supplied by Lemma~\ref{lem:movingsection} for the single
symplectic factor, transported to the total space of $\pi$, normalized so that
$q(a_t)=1$. Choose $t_i\to0$ in the open K\"ahler locus of this moving
class. Such a sequence exists because that locus is dense in $\Delta^*$.

Apply Lemma~\ref{lem:perturb} separately to $(X_{t_i},g_{t_i},a_{t_i})$ with
$\varepsilon=1/i$ and $r=i$. This gives marked $Z_i=(M,J_i)$, Gauduchon forms
$g_i$, and K\"ahler classes $a_i$ with $\rho(Z_i)=0$, $q(a_i)=1$, and
\[
d_V(a_i,a_{t_i})<1/i,\qquad
\|J_i-J_{t_i}\|_{C^i(h)}+\|g_i-g_{t_i}\|_{C^i(h)}<1/i,
\]
\[
|\mu_{Z_i,g_i}(a_i)-m(t_i)|<1/i .
\]

For later use fix a finite collection of relative holomorphic coordinate
neighbourhoods $(U_\alpha,z_\alpha)$ for the original family near $X_0$. Choose
common smaller domains $U'_\alpha\Subset U_\alpha$ on $M$, independently of $i$,
such that the $U'_\alpha$ cover $M$ and the fiber coordinate maps
$z_{\alpha,t}$ are defined near $U_\alpha$; these maps converge smoothly as
$t_i\to0$. In the local deformation of $X_{t_i}$ the coordinate theorem with
parameters extends the corresponding base point charts to smoothly varying
fiberwise holomorphic charts on the preassigned smaller domains. The additional
open conditions in Lemma~\ref{lem:perturb} allow us to require
\[
\max_\alpha\|z'_{\alpha,i}-z_{\alpha,t_i}\|_{C^i(U'_\alpha)}<1/i,
\]
with $h$ on the domains and the Euclidean norm on the coordinate spaces. The
triangle inequality then gives smooth convergence on these common domains to the
holomorphic charts of $X_0$.

In particular
\[
d_V(a_i,a_0)\le1/i+d_V(a_{t_i},a_0)\longrightarrow0 .
\]
For each fixed $r$ and every $i\ge r$ the $C^r(h)$ tensor differences are
bounded by their $C^i(h)$ differences; since the original family is smooth we
get $J_i\to J_0$ and $g_i\to g_0$ in $C^\infty$. Finally, every closed positive
$(1,1)$-current $T_i$ in $a_i$ has Gauduchon pairing $\mu_{Z_i,g_i}(a_i)$, so
\[
\int_{Z_i}T_i\wedge g_i^{d-1}=\mu_{Z_i,g_i}(a_i)\le m(t_i)+1/i\le C,
\]
which proves \eqref{eq:massbound}. The last assertion follows from the
connected component statement in Lemma~\ref{lem:perturb}.
\end{proof}

\begin{lemma}\label{lem:comparison}
Let $Z_i$ be marked irreducible compact hyperk\"ahler manifolds of dimension
$d$ with $\rho(Z_i)=0$. Suppose that $a_i,\beta_i\in\cK(Z_i)$ satisfy
\[
a_i\to a,\qquad q(a)>0,\qquad \sup_i\|\beta_i\|_V<\infty .
\]
Let $\Gamma_i$ be strictly positive $\pdbar_i$-closed real $(d-1,d-1)$-forms.
If the K\"ahler representatives $\eta_i\in a_i$ have uniformly bounded
$\Gamma_i$-mass, then so do all K\"ahler representatives of $\beta_i$.
\end{lemma}

\begin{proof}
Continuity of the fixed bilinear form on $V$ gives a constant $C_q$ with
$|q(x,y)|\le C_q\|x\|_V\|y\|_V$. Put $A=\sup_i\|a_i\|_V<\infty$ and
$B=\sup_i\|\beta_i\|_V<\infty$. Since $a_i\to a$ and $q(a)>0$ while every $a_i$
is K\"ahler, we have $c:=\inf_iq(a_i)>0$. Thus, with $K=C_qAB$, the bound
$|q(a_i,\beta_i)|\le K$ holds for every $i$. Choose $L>2K/c$. Because
$q(\beta_i)>0$, we obtain simultaneously for all $i$
\[
q(La_i-\beta_i)=L^2q(a_i)-2Lq(a_i,\beta_i)+q(\beta_i)
\ge cL^2-2KL>0,
\]
\[
q(a_i,La_i-\beta_i)=Lq(a_i)-q(a_i,\beta_i)\ge cL-K>0 .
\]
The Lorentzian signature of $q$ on $H^{1,1}(Z_i,\R)$ places $La_i-\beta_i$ in
$\cC^+_{Z_i}=\cK(Z_i)$.

For K\"ahler forms $\theta_i\in\beta_i$ and $\rho_i\in La_i-\beta_i$, all three
classes $a_i,\beta_i,La_i-\beta_i$ are K\"ahler, so the form
$\theta_i-L\eta_i+\rho_i$ is a $d$-exact real $(1,1)$-form and hence, by the
$\pdbar$-lemma on the K\"ahler manifold $Z_i$, equals $i\pdbar\varphi_i$ for a
smooth real function $\varphi_i$. Since $\Gamma_i$ is $\pdbar$-closed,
integration by parts gives
\[
\int_{Z_i}\theta_i\wedge\Gamma_i
=L\int_{Z_i}\eta_i\wedge\Gamma_i-\int_{Z_i}\rho_i\wedge\Gamma_i
\le L\int_{Z_i}\eta_i\wedge\Gamma_i,
\]
where we used that $\rho_i$ is a positive $(1,1)$-form and $\Gamma_i$ is
strictly positive, so $\rho_i\wedge\Gamma_i$ is a positive volume form and its
integral is nonnegative.
\end{proof}

\subsection{Meromorphic limits under an energy bound}

We recall the sequential compactness mechanism of
\cite[Proposition~3.3]{LS26}, including the local estimates needed here.
The hypothesis bounds the pullback of a K\"ahler form; the argument uses
this energy bound directly.

\begin{lemma}\label{lem:meromorphic}
Let $J_i\to J_0$ and $I_i\to I_0$ be smooth convergences of integrable complex
structures on compact smooth manifolds $M$ and $N$. Assume that holomorphic
source charts can be chosen to converge smoothly on common relatively compact
domains. Let $g_i\to g_0$ be Hermitian forms on $(M,J_i)$ and let
$\omega_i\to\omega_0$ be K\"ahler forms on $(N,I_i)$, including at $i=0$.
Suppose that $d=\dim_\C M\ge2$ and that holomorphic maps
$f_i\colon(M,J_i)\to(N,I_i)$ satisfy
\begin{equation}\label{eq:energy}
\sup_i\int_M f_i^*\omega_i\wedge g_i^{d-1}<\infty.
\end{equation}
Then, after passage to a subsequence, there are a proper analytic subset
$A\subset(M,J_0)$ and a meromorphic map
$f\colon(M,J_0)\dashrightarrow(N,I_0)$ such that $f_i$ converges smoothly to
$f$ on compact subsets of $M\setminus A$.
\end{lemma}

\begin{proof}
Set $T_i=f_i^*\omega_i$. These are positive closed $(1,1)$-forms for $J_i$.
The mass bound and the smooth convergence of $g_i$ give a uniform bound for
their norms as two-form measures on the fixed smooth manifold $M$. Passing to
a subsequence, we have $T_i\rightharpoonup T$. Closedness passes to the limit.
The error in testing type and positivity with $J_0$ instead of $J_i$ is bounded
by $C\|J_i-J_0\|_{C^0}$ times the mass, so $T$ is a positive closed
$(1,1)$-current on $(M,J_0)$.

We record the local estimate used to obtain convergence. In a holomorphic
source chart give the domain its Euclidean K\"ahler metric $\beta$ and put
$e_i=\tfrac12|df_i|^2$. Holomorphic maps between K\"ahler manifolds are
harmonic. The Bochner formula and the uniform bounds for the target curvature
give $\Delta e_i\ge-Ce_i^2$, with $C$ independent of $i$. Their energy and
pullback mass satisfy
\begin{equation}\label{eq:energyidentity}
e_i\,dV_\beta
=\frac{1}{(d-1)!}\,T_i\wedge\beta^{d-1}.
\end{equation}
The monotonicity formula for smooth stationary harmonic maps gives
monotonicity of
$r^{2-2d}\int_{B_r(x)}e_i\,dV_\beta$; equivalently, this follows from the
mass monotonicity formula for the positive closed current $T_i$.

These two facts imply constants $\varepsilon_*,C_*>0$, independent of $i$,
such that
\begin{equation}\label{eq:smallenergy}
r^{2-2d}\int_{B_r(x)}e_i\,dV_\beta<\varepsilon_*
\quad\Longrightarrow\quad
\sup_{B_{r/2}(x)}e_i
\le C_*r^{-2d}\int_{B_r(x)}e_i\,dV_\beta.
\end{equation}
For completeness, rescale to $r=1$. The usual weighted-maximum argument,
applied to $(\tfrac34-|z|)^2e_i(z)$, gives an a priori bound for $e_i$ on
$B_{1/2}$ when the total energy is small. Indeed, if this weighted maximum
were large, rescaling about a maximum point by $e_i^{-1/2}$ would produce a
nonnegative function with value one at the origin, bounded by four on a fixed
ball, and satisfying $\Delta v\ge-4Cv$. The scalar mean-value inequality
would give a fixed positive lower bound for its integral. Monotonicity about
the maximum point, using an outer ball of fixed radius contained in $B_1$,
would contradict sufficiently small $\int_{B_1}e_i$. Once the pointwise bound
is known, the same scalar inequality for $\Delta e_i\ge-C'e_i$ gives a linear
bound by the energy on a smaller ball. Applying this argument to fixed-radius
balls centered in $B_{1/2}$ gives \eqref{eq:smallenergy} with the displayed
radii. This derivation uses only smooth harmonic maps and uniform target
geometry.

The source charts converge smoothly, so the corresponding coordinate
expressions of $T_i$ converge weakly to that of $T$. At every radius having
zero limiting boundary mass, \eqref{eq:energyidentity} identifies the limit
of the local energies with the corresponding mass of $T$. The small-radius
limit of their normalized masses is a fixed positive dimensional multiple of
the Lelong number $\nu(T,x)$. Choose $\delta>0$ sufficiently small relative to
$\varepsilon_*$. By Siu's theorem \cite{Siu74},
\[
A=\{x\in(M,J_0):\nu(T,x)\ge\delta\}
\]
is a proper analytic subset. At every $x\notin A$ there is a ball, chosen at a
radius with zero limiting boundary mass, on which \eqref{eq:smallenergy}
applies for all sufficiently large $i$. It gives a uniform gradient bound.
Elliptic regularity for the harmonic-map equation, followed by a diagonal
argument on $M\setminus A$, gives smooth convergence to a holomorphic map
$f\colon M\setminus A\to(N,I_0)$. Lower semicontinuity gives locally finite
energy for $f$ across $A$.

It remains to obtain a meromorphic extension. Work in a relatively compact
polydisc in $(M,J_0)$. For each coordinate direction, Fubini's theorem gives
finite energy on almost every coordinate disc. A coordinate disc not contained
in $A$ meets $A$ in a discrete set. The finite-energy removable singularity
theorem \cite[Theorem~3.6]{SU81} extends the restriction of $f$ holomorphically
over each of these punctures. The complement of $A$ is a full subset of the
polydisc in the terminology of \cite{Shiffman94}. Consequently
\cite[Theorem~1]{Shiffman94}, using the Hartogs extension property of
compact K\"ahler targets \cite{Ivash92}, gives a meromorphic map on the polydisc agreeing
with $f$ almost everywhere, hence agreeing wherever both are holomorphic.
Uniqueness on overlaps glues these local extensions to the asserted
meromorphic map on $M$.
\end{proof}

\subsection{Compactness of normalized holomorphic forms}

The first mass of a Ricci-flat metric controls holomorphic two-forms in
$L^1$. This suffices for compactness of the forms; a pointwise comparison of
the Ricci-flat metrics with a background metric is unnecessary.

\begin{lemma}\label{lem:formcompact}
Let $J_i\to J_0$ be a smooth convergence of integrable complex structures on
a compact manifold $M$, and assume that holomorphic charts converge smoothly
on common relatively compact domains. Let $g_i\to g_0$ be Hermitian forms,
and let $\eta_i$ be K\"ahler forms on $(M,J_i)$ satisfying
\[
\sup_i\int_M\eta_i\wedge g_i^{d-1}<\infty.
\]
Suppose that $\sigma_i$ is a holomorphic two-form on $(M,J_i)$ and that
$\sup_M|\sigma_i|_{\eta_i}\le C$. Then a subsequence of $\sigma_i$ converges
in $C^\infty$ to a holomorphic two-form on $(M,J_0)$.
\end{lemma}

\begin{proof}
Diagonalize $\eta_i$ relative to $g_i$ at a point, with positive eigenvalues
$\lambda_1,\ldots,\lambda_d$. The induced norm on covariant two-forms gives
\[
|\sigma_i|_{g_i}
\le |\sigma_i|_{\eta_i}\max_{j<k}(\lambda_j\lambda_k)^{1/2}
\le C\,\tr_{g_i}\eta_i.
\]
Consequently
\begin{equation}\label{eq:sigmaLone}
\int_M|\sigma_i|_{g_i}\,g_i^d
\le Cd\int_M\eta_i\wedge g_i^{d-1}\le C'.
\end{equation}
In each of the converging coordinate charts, the coefficients of $\sigma_i$
are holomorphic functions and their $L^1$ norms are uniformly bounded by
\eqref{eq:sigmaLone}. The mean-value inequality and Cauchy estimates give
uniform $C^k$ bounds on every smaller chart, for every $k$. A finite covering
and a diagonal subsequence yield smooth convergence on $M$. The limiting
coefficients are holomorphic in the limiting charts.
\end{proof}

\subsection{Proof of Theorem~\ref{thm:ihs}}

\begin{proof}[Proof of Theorem~\ref{thm:ihs}]
Take $Z_i$, $a_i$ and $g_i$ from Proposition~\ref{prop:approx}, and put
$d=2m$. All classes are identified by their markings with the fixed
cohomology of $M$.

\medskip
\noindent\emph{Step 1: compactness and nondegeneracy of the symplectic forms.}

Let $\eta_i\in a_i$ be the Ricci-flat K\"ahler forms supplied by Yau's
theorem. Their Gauduchon masses satisfy
\begin{equation}\label{eq:etabound}
\int_{Z_i}\eta_i\wedge g_i^{d-1}\le C.
\end{equation}
Normalize the holomorphic symplectic forms $\sigma_i$ by
\[
q([\sigma_i],[\bar\sigma_i])=1.
\]
Write $c_Z$ for the common Fujiki constant. Since $q(a_i)=1$ and
$q(a_i,[\sigma_i])=0$, comparison of the coefficient of $st$ in the Fujiki
identity for $a_i+s[\sigma_i]+t[\bar\sigma_i]$ gives
\begin{equation}\label{eq:mixednormalization}
\int_{Z_i}\sigma_i\wedge\bar\sigma_i\wedge\eta_i^{d-2}
=\frac{c_Z}{d-1},\qquad
\int_{Z_i}\eta_i^d=c_Z.
\end{equation}
By Bochner's principle $\sigma_i$ is parallel for $\eta_i$. The pointwise
identity between
$\sigma_i\wedge\bar\sigma_i\wedge\eta_i^{d-2}$ and
$|\sigma_i|_{\eta_i}^2\eta_i^d$, with its fixed positive dimensional factor,
and \eqref{eq:mixednormalization} show that $|\sigma_i|_{\eta_i}$ is bounded
by a constant independent of $i$.

Lemma~\ref{lem:formcompact} now gives, after taking a subsequence,
\begin{equation}\label{eq:sigconv}
\sigma_i\longrightarrow\sigma_\infty\quad\text{in }C^\infty.
\end{equation}
The limit is holomorphic and closed. The polarized Fujiki identity gives
\begin{equation}\label{eq:topintegral}
\int_{Z_i}\sigma_i^m\wedge\bar\sigma_i^m
=\frac{2^m c_Z}{\binom{2m}{m}}=:C_*>0.
\end{equation}
Thus $\sigma_\infty^m$ is not identically zero. In a converging holomorphic
chart the coefficient of $\sigma_i^m$ is a nowhere vanishing holomorphic
function. Hurwitz's theorem says that its limit is either nowhere vanishing
or identically zero. The latter alternative is impossible on any chart,
because $M$ is connected and the global limiting top form is nonzero.
Therefore $\sigma_\infty^m$ is nowhere vanishing, and $\sigma_\infty$ is
everywhere nondegenerate.

Passing to cohomology in \eqref{eq:sigconv} gives
\begin{equation}\label{eq:centralperiod}
q([\sigma_\infty])=0,
\end{equation}
and
\begin{equation}\label{eq:centralpairing}
q([\sigma_\infty],[\bar\sigma_\infty])=1.
\end{equation}
In particular the limiting period belongs to $D_\Lambda$.

\medskip
\noindent\emph{Step 2: comparison with a K\"ahler family.}

Let $\phi_i\colon H^2(Z_i,\Z)\to\Lambda$ be the markings used in
Proposition~\ref{prop:approx}, and let
$\phi_0\colon H^2(X_0,\Z)\to\Lambda$ be the central topological
identification. All $(Z_i,\phi_i)$ lie in a fixed connected marked component
$\cM^0_\Lambda$. Put
\[
v_i=\phi_{i,\C}([\sigma_i]),\qquad
v_0=\phi_{0,\C}([\sigma_\infty]).
\]
Then $v_i\to v_0$ and $p_i:=\C v_i\to p_0:=\C v_0$ in $D_\Lambda$.

By the componentwise surjectivity in
Theorem~\ref{thm:globalTorelli}\textup{(ii)}, there is a marked irreducible
compact hyperk\"ahler manifold $(Y,\psi)\in\cM^0_\Lambda$ of period $p_0$.
Take a small Kuranishi family $\cY\to B$ through $Y$. Local Torelli identifies
$B$ with a neighbourhood of $p_0$ in $D_\Lambda$. Hence, after omitting
finitely many terms, there are $s_i\to0$ and marked fibers
$(Y_i,\psi_i)=(\cY_{s_i},\psi_{s_i})$ of period $p_i$.

The periods $p_i$ have Picard number zero. The singleton-fiber assertion
within $\cM^0_\Lambda$,
Theorem~\ref{thm:globalTorelli}\textup{(iii)}, therefore gives marked
biholomorphisms
\begin{equation}\label{eq:isomfi}
f_i\colon Z_i\xrightarrow{\sim}Y_i,
\qquad \phi_i\circ f_i^*=\psi_i.
\end{equation}
Choose a smooth family of K\"ahler forms $\omega_s$ on $\cY\to B$ and set
$\omega_i=\omega_{s_i}$. The classes
$\beta_i=[f_i^*\omega_i]$ are bounded in the fixed space $V$, since their
marked images are $\psi_i([\omega_i])\to\psi([\omega_0])$.
Lemma~\ref{lem:comparison}, applied with $\Gamma_i=g_i^{d-1}$, gives
\begin{equation}\label{eq:energybound}
\int_{Z_i}f_i^*\omega_i\wedge g_i^{d-1}\le C'.
\end{equation}

Choose holomorphic symplectic forms $\zeta_i$ on $Y_i$ whose marked
classes equal $v_i$. A local frame of the line of relative
holomorphic two-forms shows that $\zeta_i$ converges smoothly to a nonzero
symplectic form $\zeta$ on $Y$. The marking identity and the injectivity of
holomorphic forms into de Rham cohomology on $Z_i$ give
\begin{equation}\label{eq:pullbackzeta}
f_i^*\zeta_i=\sigma_i.
\end{equation}

Lemma~\ref{lem:meromorphic} applies to \eqref{eq:energybound}. It gives a
meromorphic map $f\colon X_0\dashrightarrow Y$ and smooth convergence on the
complement of a proper analytic subset $A\subset X_0$. Passing to the limit
in \eqref{eq:pullbackzeta} yields
\[
f^*\zeta=\sigma_\infty\quad\text{on }X_0\setminus A.
\]
The nondegeneracy of $\sigma_\infty$ makes $f$ locally biholomorphic there.
Its meromorphic graph is proper over both compact manifolds, so $f$ is
dominant and generically finite.

\medskip
\noindent\emph{Step 3: the generic degree is one.}

Resolve the graph of $f$. The image in $Y$ of the inverse image of $A$ and
of the exceptional locus over $X_0$ is a proper analytic subset. Choose a
general regular value outside this subset. If $\deg f\ge2$, it has two
distinct preimages $x_1,x_2\in X_0\setminus A$ at which $f$ is locally
biholomorphic. Choose disjoint source coordinate balls around these points
and a smaller common target ball on which both local degrees are one.
Smooth convergence of $f_i$ on smaller relatively compact source balls
preserves these degrees. Every point of a still smaller target ball then
has a preimage in each source ball for all sufficiently large $i$,
contradicting the injectivity of $f_i$. Hence $\deg f=1$, and $f$ is
bimeromorphic.

This proves $X_0\in\cC$. Ehresmann's theorem gives simple connectedness.
Bimeromorphic invariance of holomorphic forms on smooth compact complex
manifolds gives $h^{2,0}(X_0)=h^{2,0}(Y)=1$.
\end{proof}

\section{Proof of the main theorem}\label{sec:main}

\begin{proof}[Proof of Theorem~\ref{thm:main}]
We may assume that the relative dimension is positive. Take the finite \'etale
cover $p\colon\tX\to X$ from Theorem~\ref{thm:productcover}. By
Theorem~\ref{thm:product} the central fiber $\tX_0$ is the product of the central
fibers $Q_{i,0}$ of the extended factor families.

If a torus factor occurs, the extended one-forms restrict to $Q_{0,0}$ and give a
closed holomorphic coframe, since their wedge is a nowhere vanishing factor
volume form; the restriction is performed by fixing points in the other central
factors. The coframe gives a flat K\"ahler metric on $Q_{0,0}$. On its universal
cover, integration of the coframe is a local isometry to $\C^{d_0}$ and, by
completeness, a covering map; it is therefore an isomorphism. Deck transformations
preserve the pulled back coframe, so they act by translations, and compactness
gives a full lattice. Thus $Q_{0,0}$ is a complex torus.

Each irreducible symplectic factor has central fiber
bimeromorphic to an irreducible compact hyperk\"ahler manifold by
Theorem~\ref{thm:ihs}.

For a strict Calabi--Yau factor the
punctured fibers have trivial canonical bundle and $h^{0,2}=0$, so
Corollary~\ref{cor:strictCY} gives a Moishezon central fiber, which is a member
of Fujiki's class $\cC$.

The product of compact K\"ahler models of the central factors
is a compact K\"ahler model of $\tX_0$. By the K\"ahler current characterization
of class $\cC$ \cite[Theorem 0.7]{DP04}, $\tX_0$ therefore carries a K\"ahler
current. Its pushforward along the finite \'etale map $\tX_0\to X_0$ is again a
K\"ahler current: locally one sums its restrictions over the sheets of the
covering, including their smooth strictly positive lower bounds. Thus
$X_0\in\cC$.
\end{proof}

\end{document}